\documentclass[%
 aip,
 amsmath,amssymb,
 reprint,nofootinbib%
]{revtex4-1}

\usepackage{comment}
\usepackage{amsfonts}
\usepackage{tikz}
\usepackage[utf8]{inputenc}
\usepackage{graphicx}

\usepackage[caption=false]{subfig}
\usetikzlibrary{patterns}
\usepackage{amsthm}
\usepackage{float}

\usepackage{hyperref}
\DeclareMathAlphabet{\pazocal}{OMS}{zplm}{m}{n}

\usepackage{cleveref}

\usepackage{dcolumn}
\usepackage{overpic}
\usepackage{bm}

\usepackage[utf8]{inputenc}
\usepackage[T1]{fontenc}
\usepackage{mathptmx}

\usepackage{etoolbox}
\usepackage{amsmath}
\usepackage{mathtools}
\usepackage{graphicx}
\usepackage{calrsfs}

\usepackage{mathrsfs}
\usepackage{amsfonts}       

\usetikzlibrary{matrix,decorations.pathreplacing}

\pgfkeys{tikz/mymatrixenv/.style={decoration={brace},every left delimiter/.style={xshift=8pt},every right delimiter/.style={xshift=-8pt}}}
\pgfkeys{tikz/mymatrix/.style={matrix of math nodes,nodes in empty cells,left delimiter={[},right delimiter={]},inner sep=1pt,outer sep=1.5pt,column sep=2pt,row sep=2pt,nodes={minimum width=20pt,minimum height=10pt,anchor=center,inner sep=0pt,outer sep=0pt}}}
\pgfkeys{tikz/mymatrixbrace/.style={decorate,thick}}
\newcommand*\mymatrixbraceright[4][m]{
    \draw[mymatrixbrace] (#1.west|-#1-#3-1.south west) -- node[left=2pt] {#4} (#1.west|-#1-#2-1.north west);
}
\newcommand{\norm}[1]{\left\lVert#1\right\rVert}

\DeclareMathOperator*{\argmin}{arg\,min}
\newcommand{\cM}{\pazocal{M}}

\makeatletter
\def\@email#1#2{%
 \endgroup
 \patchcmd{\titleblock@produce}
  {\frontmatter@RRAPformat}
  {\frontmatter@RRAPformat{\produce@RRAP{*#1\href{mailto:#2}{#2}}}\frontmatter@RRAPformat}
  {}{}
}%
\makeatother

\let\svthefootnote\thefootnote
\newcommand\freefootnote[1]{%
  \let\thefootnote\relax%
  \footnotetext{#1}%
  \let\thefootnote\svthefootnote%
}
\begin{document}

\title[Learning Dynamics on Invariant Measures using PDE-Constrained Optimization]{Learning Dynamics on Invariant Measures \\Using PDE-Constrained Optimization}
\author{Jonah Botvinick-Greenhouse}
 \altaffiliation[]{Corresponding author: \texttt{jrb482@cornell.edu}}
\affiliation{Center for Applied Mathematics,
Cornell University,
Ithaca, NY 14850}

\author{Robert Martin}%
\affiliation{ 
DEVCOM Army Research Laboratory,
Research Triangle Park,
Durham, NC 27709
}%
\author{Yunan Yang}

\affiliation{%
Institute for Theoretical Studies,
ETH Z\"urich,
Z\"urich, Switzerland 8092
}%

\date{\today}

\begin{abstract}
We extend the methodology in [Yang et al., 2023] to learn autonomous continuous-time dynamical systems from invariant measures. The highlight of our approach is to reformulate the inverse problem of learning ODEs or SDEs from data as a PDE-constrained optimization problem. This shift in perspective allows us to learn from slowly sampled inference trajectories and perform uncertainty quantification for the forecasted dynamics. Our approach also yields a forward model with better stability than direct trajectory simulation in certain situations. We present numerical results for the Van der Pol oscillator and the Lorenz-63 system, together with real-world applications to Hall-effect thruster dynamics and temperature prediction, to demonstrate the effectiveness of the proposed approach.
\end{abstract}

\maketitle

\begin{quotation}
Data-driven models have proven to be instrumental across numerous scientific disciplines for their ability to predict and control the behavior of complex physical systems.\cite{MONTANS2019845} Popular approaches for modeling dynamical trajectories typically adopt a Lagrangian perspective and seek a pointwise matching with either the observed data or its approximate state derivatives. When the observed data has a poor temporal resolution and the state derivatives are difficult to approximate, these approaches may struggle. Such difficulties are further exaggerated when measurements are contaminated with noise, and the system in question exhibits sensitive dependence on initial conditions. In this paper, we propose an alternative approach that can circumvent some of these challenges by treating global statistics of the observed dynamics as the inference data.
\end{quotation}
\section{Introduction}\label{sec:Intro}

\begin{figure*}[!ht]
    \centering
    \includegraphics[width = \textwidth]{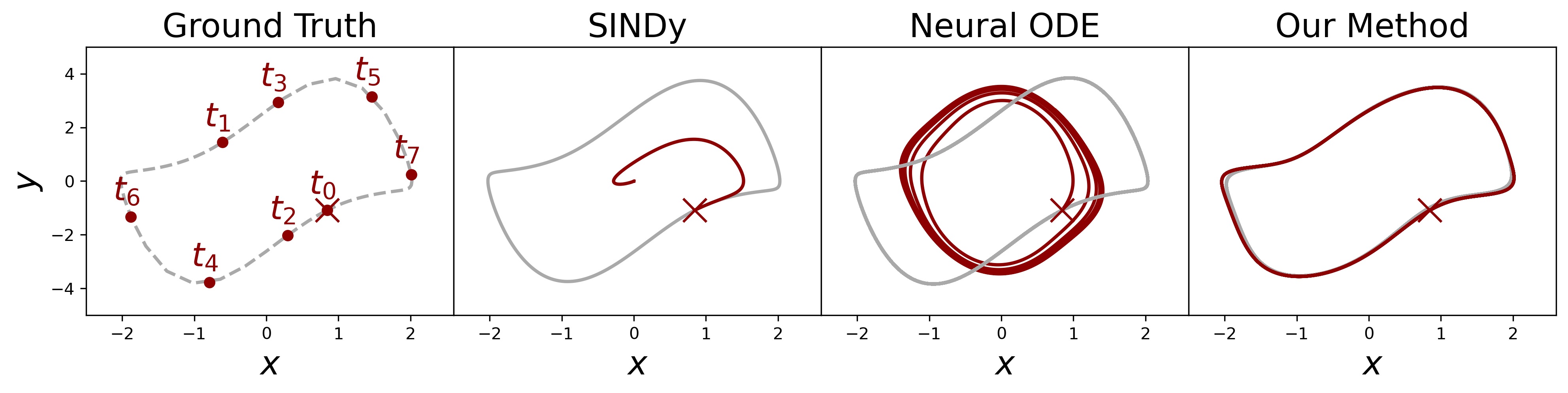}
    \caption{Comparison with the SINDy \cite{brunton2016discovering,desilva2020,Kaptanoglu2022} and the Neural ODE \cite{DBLP:journals/corr/abs-1806-07366} frameworks for learning slowly sampled dynamics. The left panel shows the original dynamics (in grey) and the first eight points of a slowly sampled trajectory (in red). While the SINDy and Neural ODE frameworks can learn the quickly sampled dynamics (model output in grey), both struggle to learn from the slowly sampled trajectory (model output in red). On the other hand, our framework can learn from both quickly and slowly sampled dynamics (the red and grey model outputs coincide). Additional experiment details are provided in \Cref{subsec:vanderpol}.}
    \label{fig:initcompare}
\end{figure*}

Differential equations are typically used to model trajectory data originating from physical systems. Common techniques for fitting differential equations to trajectory data include the shooting methods, \cite{baake1992fitting,michalik2009incremental} neural differential equations, \cite{DBLP:journals/corr/abs-1806-07366, DBLP:journals/corr/abs-1905-10403, DBLP:journals/corr/abs-2009-03288} and SINDy.\cite{brunton2016discovering,fasel2022ensemble} Methods based on the Kalman filter are effective in data assimilation and in estimating unknown states and parameters of a system as it evolves in time.~\cite{van2001square,evensen2003ensemble,simon2006optimal,harlim2014ensemble,hamilton2016ensemble,schillings2017analysis} When used to identify model parameters, these approaches fall under the broad category of system identification.

These approaches all adopt a Lagrangian perspective and directly fit the modeled trajectories or their state derivatives to the observed measurements. While these techniques have seen great success in modeling complex dynamics, their application is generally limited to inference trajectories sampled at a relatively high frequency. When the inference trajectory is sampled slowly, or in the worst-case scenario when measurement times are unknown, these approaches may not be applicable. For example, see Figure \ref{fig:initcompare} in which we investigate the use of SINDy\cite{brunton2016discovering} and a Neural ODE\cite{DBLP:journals/corr/abs-1806-07366} for modeling the dynamics of a slowly sampled limit cycle. 

 There are at least three sources of instabilities when directly using the trajectory data to perform velocity reconstruction. First, for certain chaotic dynamical systems, a small perturbation in the initial condition can lead to a large deviation in the trajectory at a later time, which cannot be differentiated from inaccurate dynamics by looking at the data alone. Second, the estimation of the particle velocity suffers from slowly sampled trajectory data, which directly affects the reconstruction of the dynamics, as shown in Figure~\ref{fig:initcompare}. Third, the measurement (extrinsic) noise and the model intrinsic noise both change the state location. The small noise pollution can be amplified more in the velocity estimation using the divided difference with a small time step. All three factors share the nature that a small perturbation in the trajectory data leads to a large deviation in the estimated velocity/learned dynamics. 

In contrast with the Lagrangian approach to modeling dynamics, our method builds on an Eulerian perspective\cite{greve2019data,yang2021optimal} in which velocity models are constructed to yield the same asymptotic statistics as the observed measurements. This approach converts what is traditionally regarded as an ordinary differential equation (ODE) or stochastic differential equation (SDE) modeling problem into a partial differential equation (PDE)-constrained optimization problem. The motivation of our method is that, in certain situations, the PDE forward model yields better stability in solving the inverse problem than direct trajectory forward simulation based upon an ODE or SDE. Importantly, our method does not rely on prior knowledge of sampling times and can thus be used to learn the dynamics from slowly sampled trajectories.  

There are two important differences between the line of work using Kalman filters and our proposed method. First, a Kalman filter is a particular case of the Bayes filter using the Bayes theorem, while our reconstruction follows a deterministic inverse problem (PDE-constrained optimization). Second, time is a crucial element in designing a Kalman filter, while in our approach, we use the invariant measure and a time-independent PDE surrogate model.  Once the flow has been inferred, we can also perform uncertainty quantification for the forecasted dynamics, building towards extending grid-based Bayesian estimation of nonlinear low-dimensional systems\cite{bewley2012efficient} to slowly sampled unknown systems with nontrivial invariant measures.

More specifically, instead of directly treating the noisy observations $\{\tilde{\mathbf{x}} (t_i)\}_{i=1}^n$  from one single trajectory of an autonomous flow $\dot{\mathbf{x}} = v^*(\mathbf{x})$ as inference data, we consider the occupation measure $\rho^*$ generated by a single trajectory, where for each measurable set $B$,
\begin{equation}\label{E1}
\rho^*(B) := \frac{1}{n}\sum_{i=1}^n \chi_B\left({\tilde{\mathbf{x}} (t_i)}\right), \hspace{.6cm} \chi_B(x) = \begin{cases} 1, & x\in B, \\
0 ,& x\not\in B.
\end{cases}
\end{equation}
When the occupation measures generated by a nontrivial (see \Cref{subsec:PM}) set of initial conditions all weakly converge to the same invariant measure, the limiting measure is said to be physical.\cite{SRB} In this work, we consider the class of autonomous systems for which the occupation measure of Lebesgue-almost all initial conditions converges to a unique physical measure. Notably, this encompasses chaotic attractors such as the Lorenz-63 system.\cite{luzzatto2005lorenz,TUCKER19991197} This assumption guarantees the uniqueness of the invariant measure for the dynamical system under study. If we relax it and allow the existence of multiple invariant measures, further treatment of the PDE forward model is needed. For instance, the fact that different invariant measures are mutually singular as well as information on the initial condition, among other considerations, is necessary to guarantee that the steady-state solution picked up by the PDE model matches the observed invariant measure. We remark that the definition of a physical measure demonstrates its robustness to perturbations with respect to initial conditions. 

Going forward, we write $v = v(\theta) = v(\mathbf{x};\theta)$ to denote the dependence of the reconstructed velocity fields on a set of parameters $\theta\in\Theta$ where $\Theta\subset \mathbb{R}^m$ is the admissible set of all parameter values. The concrete form of $\theta$ depends on the hypothesis space of $v$, which will be discussed in~\Cref{sec:velocity_parameter}. The task is now to find the best-parameterized model $v(\mathbf{x};\theta)$ approximating the true velocity $v^*$ through the optimization
\begin{equation}\label{E2}
    \inf\limits _{\theta \in \Theta} \pazocal{J}(\theta), \hspace{1cm} \pazocal{J}(\theta):= \pazocal{D}(\rho_{\varepsilon}(v(\theta)),\rho^*).
\end{equation} 
The formulation~\eqref{E2} represents an inverse data-matching problem, in which $\pazocal{D}$ denotes a metric or divergence on the space of probability measures and $\rho_{\varepsilon}(v(\theta))$ is a regularized approximation to the physical measure of the dynamical system, given some regularization parameter $\varepsilon > 0$ and the current velocity $v(\theta)$. That is, $v(\theta) \mapsto \rho_{\varepsilon}(v(\theta))$ is our new forward model.

Although one could approximate $\rho(v(\theta))$ by numerically integrating a trajectory and binning the observed states to a histogram,\cite{greve2019data} this approach does not permit simple differentiation of the resulting measure with respect to the parameters $\theta$. When the size of $\theta$, i.e., $m$, is large, it is practical to use gradient-based optimization methods for solving the optimization problem~\eqref{E2}, and one has to compute the essential gradient $\partial_{\theta} \pazocal{J}$. In Ref.~\onlinecite{yang2021optimal}, this was handled by viewing $\rho_{\varepsilon}(v(\theta))$ as the dominant eigenvector of a regularized Markov matrix originating from an upwind finite volume discretization of the continuity equation. The derivative $\partial_{\theta} \pazocal{J}$ was then seamlessly computed via the adjoint-state method.\cite{yang2021optimal} The computation time of the adjoint-state method is independent of the size of $\theta$, making the framework presented in Ref.~\onlinecite{yang2021optimal} well-suited for large-scale computational inverse problems.

In this work, we build upon the framework proposed in Ref.~\onlinecite{yang2021optimal} and study invariant measure-based velocity learning with a large-scale parameter space applied to real data. There are three essential new contributions: 
\begin{enumerate}
    \item We consider the Fokker--Planck equation as the partial differential equation (PDE) forward model for $\rho_{\varepsilon}(v(\theta))$, rather than the continuity equation. This is motivated by the Fokker--Planck equation's greater modeling capacity. Indeed, the Fokker--Planck equation reduces to the continuity equation when its diffusion term is zero, and it can fit intrinsic noise present in trajectories which reduces over-fitting the parameterized velocity $v(\theta)$. Moreover, the Fokker--Planck equation can be seen as an alternative to the teleportation regularization used for the continuity equation in Ref.~\onlinecite{yang2021optimal} in order to guarantee the uniqueness of the computed stationary solution $\rho_{\varepsilon}(v(\theta))$. 
    \item In contrast to only learning three coefficients as done in Ref.~\onlinecite{yang2021optimal}, we parameterize the velocity $v(\theta)$ using piecewise polynomial, global polynomial, and neural network discretizations, which can all yield large parameter spaces with thousands of dimensions. We compare the reconstructed velocity in each case and further discuss how the choice of parameterization affects the inverse problem's well-posedness and the reconstructed velocity's regularity. We also consider various metrics/divergences as the choice of the objective function.
    \item We investigate velocity learning in time-delay coordinates, which can characterize the full dynamics from partial state measurements alone.\cite{Takens1981DetectingSA} After performing the optimization \eqref{E2}, we evolve the learned Fokker--Planck equation forward in time to quantify the uncertainty in predictions of future dynamics. Based on this framework, we demonstrate that forecasts incur larger uncertainties when the embedding dimension is not sufficiently high. It is worth noting that there is no analytic form for the velocity in time-delay coordinates, even for well-studied dynamical systems. We also stress that our proposed approach permits larger-scale modeling of time-delayed dynamics than the approach considered in Ref.~\onlinecite{greve2019data}, due to the use of the adjoint-state method when solving the PDE-constrained optimization.
\end{enumerate}

The rest of the paper is organized as follows. In~\Cref{sec:background}, we review essential background on dynamical systems, invariant measures, the Fokker--Planck equation, and time-delay coordinates. In~\Cref{sec:forward}, we introduce the forward surrogate model $\rho_\varepsilon(v(\theta))$ and analyze its modeling errors. In~\Cref{sec:gradientcalc}, we present an efficient gradient calculation for the objective function $\pazocal{J}(\theta)$ by treating~\eqref{E2} as a PDE-constrained optimization problem and utilizing the adjoint-state method. We then adapt the gradient calculation to various velocity parameterizations, including neural network discretizations in which the gradient is computed along with the backpropagation technique.\cite{lecun1989backpropagation} 

Finally, in~\Cref{sec:numerics}, we present velocity reconstructions for the Van der Pol oscillator and the Lorenz-63 system. We also model dynamics in time-delay coordinates based on real-world data from a Hall-effect thruster and actual temperature recordings. We perform uncertainty quantification on the last two real-data examples. Conclusions follow in~\Cref{sec:conclusions}.

\section{Background}\label{sec:background}
This section reviews the essential background on invariant measures, stochastic dynamics, the Fokker--Planck equation, and time-delay coordinates. We also review the Eulerian approach for parameter identification proposed in Refs.~\onlinecite{greve2019data,yang2021optimal}, as well as past work on the discrete inverse Frobenius--Perron problem.\cite{IFP1} 

\subsection{Physical Measures}\label{subsec:PM}
Physical measures characterize the long-term statistical behavior of a significant collection of dynamical trajectories. When a dynamical system is chaotic and exhibits sensitive dependence on initial conditions, the existence of a physical measure unifies the statistical properties of trajectories that are pointwise dissimilar. While ergodic measures also describe the long-term statistical behavior of dynamical trajectories, they may have very small support or even be singular. On the other hand, when a dynamical system admits a physical measure, it holds that the trajectories corresponding to a positive Lebesgue measure subset of initial conditions will all share the same statistical behavior. 
We will now formalize these ideas in the language of ergodic theory. For a more thorough treatment of the topic, we refer to Refs.~\onlinecite{cowieson_young_2005, Froyland96estimatingphysical,SRB,katok_hasselblatt_1995}.

While we will review the theory of physical measures in the context of discrete-time dynamical systems, our applications will consider dynamics given by a time-$\Delta t$ flow map for some $\Delta t>0$. Following Ref.~\onlinecite{SRB}, we assume that $\cM$ is a compact Riemannian manifold and that $T:\cM\to \cM$ is a diffeomorphism. A probability measure $\mu$ is said to be invariant with respect to the map $T$ if $\mu(T^{-1}(B)) = \mu(B)$ for all $B\in\mathscr{B}$, where $\mathscr{B}$ denotes the Borel $\sigma$-algebra (see Ref.~\onlinecite[Definition 2.1]{EW}). Hereafter, we will assume that $\mu$ is an invariant measure. A point $x\in \cM$ is said to be generic (see Ref.~\onlinecite[Section 2.2]{SRB}) if for all $g\in C(\cM)$, it holds that 
\begin{equation}\label{eq:generic}
    \lim_{N\to \infty}\frac{1}{N}\sum_{k=0}^{N-1}g(T^k(x)) = \int_{\pazocal{M}} g\, d\mu.
\end{equation}
The left-hand side of \eqref{eq:generic} is known as the time-average of a function $g\in C(\cM)$ whereas the right-hand side of \eqref{eq:generic} is known as the space average. It follows Birkhoff's pointwise ergodic theorem (see Ref.~\onlinecite[Theorem 2.30]{EW}) that the time-average of any $g\in C(\cM)$ necessarily exists on a set of full $\mu$-measure. To formally discuss the statistical properties of dynamical trajectories, we now define the $N$-step occupation measure given the initial condition $x\in \cM$ as
\begin{equation}\label{eq:occupation}
    \mu_{x,N}(B):= \frac{1}{N}\sum_{k=0}^{N-1}\chi_B(T^k(x)),\hspace{1cm} \forall B\in\mathscr{B}.
\end{equation}
The condition that a point $x\in \cM$ is generic is equivalent to the condition 
\begin{equation}\label{eq:weakstar}
    \lim_{N\to\infty}\mu_{x,N}=\mu\,,
\end{equation} where convergence takes place in the weak-* topology (see Ref.~\onlinecite[Definition 4.19]{EW}). Since the quantity $\mu_{x,N}(B)$ approximates the average amount of time for which the orbit $\{T^k(x)\}_{k=0}^{\infty}$ initiated at $x\in\cM$ resides in a measurable set $B\in\mathscr{B}$, this convergence indicates that the collection of generic points all share the same asymptotic statistical behavior. When the measure $\mu$ is ergodic (see Ref.~\onlinecite[Definition 4.19]{EW}), it holds that $\mu$-almost every $x\in \cM$ is a generic point (see Ref.~\onlinecite[Corollary 4.20]{EW}). However, if $\mu$ is an ergodic measure that is singular with respect to the Lebesgue measure, the resulting collection of generic points may be physically insignificant and difficult to observe computationally. Motivated by this perspective, an invariant measure $\mu$ is said to be physical if there exists a collection of generic points with positive Lebesgue measure (see Ref.~\onlinecite[Definition 2.3]{SRB}).

We will next discuss the ways in which a physical invariant measure $\mu$ can be computationally approximated. If one collects the measurements $\{T^k(x)\}_{k=1}^N$, the weak-* convergence in \eqref{eq:weakstar} suggests that the physical measure $\mu$ will describe the statistics of our measurements provided that $N$ is sufficiently large. Motivated by this perspective, we can discretize the domain $\cM$ and directly compute the occupation measure \eqref{eq:occupation} for each cell in the discretization to approximate the physical measure. This procedure has been previously used to approximate physical measures.\cite{Allawala_2016,greve2019data,yang2021optimal} Other approaches have been proposed to compute the invariant measure as the stationary vector of the finite-dimensional approximation of the continuous Frobenius--Perron operator,\cite{lasota1998chaos} including Ulam’s method,\cite{Froyland96estimatingphysical} and Galerkin-type methods.\cite{10.1007/978-3-642-56589-2_7,Sch_tte_2016} More precisely, these discretizations are used to construct a Markov matrix that represents a random dynamical system approximating the deterministic map $T:\cM\to \cM$. An invariant measure for the discrete approximation is then recovered as a stationary vector of the resulting Markov matrix. As the discretization is refined, certain assumptions guarantee that the desired physical measure will be recovered in the weak-* limit (see Ref. \onlinecite[Theorem 4.14]{10.1007/978-3-642-56589-2_7}). 

\subsection{Stochastic Dynamics and the Fokker--Planck Equation}
Consider an It\^o stochastic differential equation (SDE) of the form 
\begin{equation}\label{ES}
 dX_t = v(X_t) dt+\sigma(X_t) dW_t, \hspace{1cm} X_0 = x.   
\end{equation}
Above, $W_t$ is a Brownian motion, $v$ is the velocity, and $\sigma$ determines the diffusion matrix $\Sigma (\mathbf{x}) = \frac{1}{2} \sigma(\mathbf{x}) \sigma(\mathbf{x})^\top.$ For simplicity, we will consider the case of a constant diffusion. Similar to the deterministic setting, there are analogous notions of invariant measures, ergodicity, and physical measures in the stochastic setting. \cite{Blumenthal_2019,Hong2019} One may use the Euler--Maruyama method to obtain the numerical solution to~\eqref{ES} on the time interval $[0,T]$, which assigns
$$X_{j+1} = X_j + v(X_j) \Delta t +\sigma(X_j) \xi_j \sqrt{\Delta t},$$
where $\{\xi_j\}$ are independently and identically distributed (i.i.d.) from $\pazocal{N}(0,I)$, the standard normal distribution on $\mathbb{R}^d$, $\Delta t:=T/N,$ and $j \in \{0,\dots,N-1\}$.

The Fokker--Planck equation provides a PDE description of the probability density $\rho(\mathbf{x},t)$ of the random variable $X_{t}$. The density evolves as (see Ref.~\onlinecite[Page 88]{SP})
\begin{equation}\label{FP1}
    \frac{\partial \rho(\mathbf{x},t)}{\partial t} = -\nabla \cdot (\rho(\mathbf{x},t) v(\mathbf{x})) + \nabla \cdot \Big(\nabla \cdot(\Sigma(\mathbf{x}) \rho(\mathbf{x},t))\Big).
\end{equation}
By assuming a constant diffusion, we may write $\Sigma(\mathbf{x}) = D I,$ where $I$ denotes the identity and $D >0$ is a constant representing the scale of the diffusion. Equation \eqref{FP1} can then be simplified to read 
\begin{equation}\label{FP2}
    \frac{\partial \rho(\mathbf{x},t)}{\partial t} = -\nabla \cdot (\rho(\mathbf{x},t) v(\mathbf{x})) + D  \nabla^2 \rho(\mathbf{x},t).
\end{equation}
We leave the study of a non-constant or anisotropic diffusion for later work. 
We remark that if $D = 0$, \eqref{FP2} reduces to the so-called continuity equation, which instead models the probability flow of the ODE given by $\dot{\mathbf{x}} = v(\mathbf{x}).$ Under certain conditions,\cite{huang2015steady} the steady-state solution $\rho(\mathbf{x})$ of~\eqref{FP2} exists and satisfies
\begin{equation}\label{eq:FPE}
    \nabla \cdot (\rho(\mathbf{x}) v(\mathbf{x})) = D \nabla^2 \rho(\mathbf{x}).
\end{equation}
Since~\eqref{eq:FPE} describes a limiting distribution $\lim_{t\to \infty} \rho( \mathbf{x},t),$ it has been previously used to provide approximations of invariant measures for stochastically-forced dynamical systems.\cite{Allawala_2016} In Ref.~\onlinecite{chen2021solving}, an SDE learning problem was studied using~\eqref{FP1} as the modelling equation with different data assumptions.

\subsection{Delay Coordinates and Takens' Theorem}
The technique of time-delay embedding is a popular approach for reconstructing chaotic dynamical systems from limited observations.\cite{Brunton_2017,greve2019data,kirtland2022unstructured, Sugihara2012} The procedure involves embedding time series measurements $\psi(t)=\psi(\mathbf{x}(t))$ of the state $\mathbf{x}(t)$ into $d$-dimensional Euclidean space by considering the vector of time-lagged observations $$\Psi_{d,\tau}(t)=(\psi(t),\psi(t-\tau),\dots, \psi(t-(d-1)\tau)),$$ for some $\tau > 0$. Takens' theorem \cite{Takens1981DetectingSA} provides suitable assumptions under which $\Psi_{d,\tau}(t)$ and $\mathbf{x}(t)$ are related via a diffeomorphism, implying that the time-lagged vector of partial observations $\Psi_{d,\tau}(t)$ is sufficient for reconstructing the full state $\mathbf{x}(t).$ Notably, the embedding dimension provided in Ref.~\onlinecite{Takens1981DetectingSA} is $d=2m+1$ where $m$ is the dimension of a compact manifold $\pazocal{M}$ on which the flow map $f_t$ for the original dynamics is defined. In cases when trajectories are attracted to a compact subset $A$ with box-counting dimension (see Ref.~\onlinecite[Page 586]{Sauer1991}) $d_A$ strictly less than $m$, it turns out that lower-dimensional embeddings can be obtained. 

When a time-series projection $\psi(t)$ of an unknown system $\dot{\mathbf{x}} = v(\mathbf{x})$ is observed, one can try to numerically determine a suitable embedding dimension $d$ and time delay $\tau$; see for example, Refs.~\onlinecite{10.1115/1.4036814,Hongguang2006SelectionOE, NNdelay,  Wallot2018CalculationOA}. Choosing a proper embedding dimension and time delay is important for obtaining a reliable surrogate model of the original dynamics in time-delayed coordinates. Notably, in \Cref{subsec:HET}, we demonstrate that models for the velocity in time-delayed coordinates can incur excess uncertainties when the embedding dimension is not sufficiently large. 

\subsection{Prior Work on Learning Dynamics from Invariant Measures}\label{subsec:prior}
For chaotic systems, trajectories are sensitive to initial conditions and estimation parameters. Sometimes, the approximate reference velocity field $\{\hat{v}(\mathbf{x}(t_i))\}$ cannot be accurately estimated from a trajectory $\{\mathbf{x}(t_i)\}$ due to the lack of observational data, slow sampling, discontinuous or inconsistent time trajectories, and noisy measurements. To tackle such difficulties, instead of working with the Lagrangian trajectories, Refs.~\onlinecite{greve2019data,yang2021optimal} propose an Eulerian approach by treating the occupation measure~\eqref{eq:occupation} as the data. When enough samples are available, the occupation measure can be treated as an approximation to the invariant measure; see~\Cref{subsec:PM}. Finding the optimal parameter $\theta$ is then translated into the optimization problem~\eqref{E2}. The reference measure $\rho^*$ is the occupation measure converted from the observed trajectories $\{\hat{\mathbf{x}}(t_i)\}$; see~\eqref{eq:occupation}. In Ref.~\onlinecite{greve2019data}, the approximated synthetic $\rho_\varepsilon(v(\theta))$ is generated by first simulating the synthetic  trajectories  $\{\mathbf{x}(t_i;\theta)\}$ based on the dynamical system and then computing its histogram following~\eqref{eq:occupation}. Since this approach requires lengthy trajectory simulation, each evaluation of $\rho_\varepsilon(v(\theta))$ for a given $\theta$ is relatively costly. Moreover, it is difficult to compute the derivative of $\rho_\varepsilon(v(\theta))$  with respect to $\theta$ due to the histogram approximation of nonlinear trajectories. As an improvement to the original idea in Ref.~\onlinecite{greve2019data}, Ref.~\onlinecite{yang2021optimal} proposes a surrogate model to approximate $\rho_\varepsilon(v(\theta))$ that is differentiable in $\theta$ and sometimes faster to compute. The key idea is to solve for $\rho_\varepsilon(v(\theta))$ as the distributional steady-state solution to the continuity equation (i.e., \eqref{eq:FPE} with $D = 0$) using a finite volume upwind scheme together with the teleportation regularization. The gradient of the objective function $\pazocal{J}$ in~\eqref{E2} with respect to the parameter $\theta$ can be efficiently computed based on the adjoint-state method (see Ref.~\onlinecite[Sec.~5]{yang2021optimal}). The problem of learning an SDE from an invariant measure is also studied in Ref.~\onlinecite{chen2023detecting}, which uses a deep learning framework to invert the drift and diffusion terms.

The task of learning a dynamical system from an invariant measure has also been studied in the discrete-time setting under the inverse Frobenius--Perron problem.\cite{IFP1, IFP2, IFP3, IFP4} The Frobenius--Perron operator, also known as the transfer operator, characterizes the time evolution of an initial measure $\mu_0$ according to some prespecified dynamical system. Given a probability measure $\mu$, the inverse Frobenius--Perron problem seeks to construct a dynamical system for which $\mu$ is a fixed point of the associated transfer operator. The most widely studied case involves recovering an ergodic map $T$ on $[0,1]$ for which a prescribed absolutely continuous measure is the unique fixed point of the discrete transfer operator. In this particular setting, various approaches such as topological conjugation \cite{grossmann1977invariant} and matrix methods \cite{NIE2018248} have been introduced to solve the inverse problem. The multivariate inverse Frobenius--Perron problem was also studied in Ref.~\onlinecite{Fox_2021}, where ergodic maps were constructed to adhere to the statistics of two-dimensional densities. Moreover, due to inherent non-uniqueness in the inverse problem, recent approaches further restrict the solution space of the discrete ergodic maps to those with a prescribed power spectrum. \cite{mcdonald2020novel} To the best of our knowledge, Refs.~\onlinecite{greve2019data,yang2021optimal,chen2023detecting}, and our contributions here are the first works that numerically solve the inverse Frobenius--Perron problem in the continuous-time setting. Notably, we do not assume that $\mu$ is absolutely continuous, as we use a finite-volume discretization to approximate the Frobenius--Perron operator.

\section{The Forward Model and Modeling Errors}\label{sec:forward}
A central contribution of this work is to consider a different regularized forward model than the one in Ref.~\onlinecite{yang2021optimal}, especially for trajectory measurements containing intrinsic noise, which can be interpreted as sample paths of stochastic dynamical systems~\eqref{ES}. In those cases, the Fokker--Planck equation~\eqref{FP1} is a better candidate as the PDE surrogate model, as it contains a diffusion term that can fit noise present in the data. Based on the relationship between~\eqref{ES} and~\eqref{FP1}, one can learn both the velocity field $v(\mathbf{x})$ and the diffusion tensor $\Sigma(\mathbf{x})$ in the optimization framework~\eqref{E2}. For simplicity, we only consider a fixed diffusion constant and leave the investigation of multi-parameter inversion to future work.

We will use~\eqref{eq:FPE} as the forward model to fit invariant measures generated by trajectories with intrinsic noise. While the diffusion term allows the model to fit the intrinsic noise and prevent over-fitting the noise into the target velocity component, it also controls the scaling of the reconstructed velocity $v(\mathbf{x};\theta)$. Indeed, when $D=0$ and $\Tilde{v}(\mathbf{x}) = a\,v(\mathbf{x})$, we have $\nabla\cdot(\rho(\mathbf{x}) \Tilde{v}(\mathbf{x} ) ) =0$ as long as $\nabla\cdot(\rho(\mathbf{x}) v(\mathbf{x} ) ) =0$, for any $a > 0$. However, for most cases,  $\widetilde{v}$ and $v$ will not solve the stationary Fokker-Planck equation \eqref{eq:FPE} for $D > 0$.

\subsection{Finite Volume Discretization}\label{sec:FVM}
We assume that our system evolves on the $d$-dimensional rectangular state space 
$$\Omega= [a_1,b_1]  \times \dots \times [a_d,b_d]  \subset \mathbb{R}^d,$$ 
with a spatially dependent velocity $v:\Omega \to \mathbb{R}^d$. We define $n_i \in \mathbb{Z}^+$, $1\leq i \leq d$, to be the number of equally-spaced points along the $i$-th spatial dimension at which we wish to approximate the solution of \eqref{FP2}, as well as the mesh spacing $$\Delta x_i:= \frac{b_i-a_i}{n_i-1}.$$ We are thus interested in obtaining a solution to the forward problem at points of the form $$x_{k_1,\dots,k_d}:= ( a_1 + k_1\Delta x_1, \dots, a_d+k_d \Delta x_d),$$
where $k_i \in \{1,\dots, n_i\}$. We will index our coordinates using column-major order and write $x_{k_1,\dots,k_d} = x_j$ where 
\begin{equation}\label{eq:Si}
j =k_1+ \sum_{i=2}^d (k_i-1)S_i\,, \hspace{.5cm}  S_{i} := \prod_{j = 1}^{i-1}n_j \,.
\end{equation}
We will regard $x_j$ as the center of the cell $C_j$ where 
$$
C_j = \prod_{i=1}^d \left[a_i +\left( k_i -\frac{1}{2} \right)\Delta x_i  , \,  a_i + \left(k_i + \frac{1}{2} \right)\Delta x_i\right).
$$   
Following the approach in Ref.~\onlinecite{bewley2012efficient}, we implement a first-order upwind finite volume discretization of the continuity equation, adding a diffusion term using the central difference scheme and enforcing a zero-flux boundary condition.\cite{leveque2002finite} This allows us to obtain an explicit time-evolution of the probability vector $\rho= \begin{bmatrix}\rho_1 & \rho_2 & \dots & \rho_N
\end{bmatrix}^{\top}\in\mathbb{R}^N$, where $N=\prod_{i=1}^d n_i$. While $\rho$ is a discrete probability measure over the cells $C_j$, it also corresponds to a piecewise-constant probability density function on $\Omega$. 

With an abuse of notation, we will refer to both the piecewise-constant density and the discrete probability measure as $\rho$.  We discretize the time domain with a step size $\Delta t$. Based on~\eqref{FP2}, the probability vector at the $l$-th time step evolves as 
$$
\rho^{(l+1)} = \rho^{(l)}+K \rho^{(l)}, \hspace{1cm} K = \sum_{i=1}^d \frac{\Delta t}{\Delta x}K_{i} ,
$$
where each $K_{i}$ is a tridiagonal matrix of the form
\begin{widetext}

\begin{equation*} 
K_{i} :=\begin{tikzpicture}[baseline={-0.5ex},mymatrixenv]
\matrix [mymatrix,inner sep=4pt](m)
{
\ddots  \\
& -v^{i,-}_{j-1} +\displaystyle \frac{D}{\Delta x_i }\\ 
\ddots& \vdots\\ & v^{i, -}_{j-1}-w^{i,+}_{j-1} -\displaystyle \frac{2D}{\Delta x_i }&  -v^{i,-}_{j} +\displaystyle \frac{D}{\Delta x_i } \\ 
\ddots & \vdots &\vdots \\ 
& w^{i,+}_{j-1}+\displaystyle \frac{D}{\Delta x_i }& v^{i,-}_{j} - w^{i,+}_{j} -\displaystyle \frac{2D}{\Delta x_i }& -v^{i,-}_{j+1}+\displaystyle \frac{D}{\Delta x_i }\\
& &  \vdots & \vdots & \ddots\\
& & w^{i,+}_{j}+\displaystyle \frac{D}{\Delta x_i } & v^{i,-}_{j+1}-w^{i,+}_{j+1}-2\displaystyle \frac{D}{\Delta x_i } \\
& & & \vdots & \ddots \\
& & &  w^{i,+}_{j+1}+\displaystyle \frac{D}{\Delta x_i } \\
& & & & \ddots \\
};
\mymatrixbraceright{1}{3}{$S_{i}$}
\end{tikzpicture}\in\mathbb{R}^{N\times N}.
\end{equation*}
\end{widetext}

Above, we have defined for each $j\in \{1,\dots, N\}$ the upwind velocities
\begin{align*}
v_{j}^{i,-} := \min \big\{0,v_j^i\big\},& \hspace{1cm}v_{j}^{i,+} := \max \big\{0,v_j^i\big\},\\
w_{j}^{i,-} := \min \big\{0,w_j^i\big\},& \hspace{1cm}w_{j}^{i,+} := \max \big\{0,w_j^i\big\},
\end{align*}
where  $ v_j^i := v\big(x_{j} - \mathbf{e}_i\Delta x_i/2\big)\cdot \mathbf{e}_i$ and $w_j^i :=v\big(x_{j} + \mathbf{e}_i\Delta x_i/2\big)\cdot \mathbf{e}_i$ denote the $i$-th components of the velocity field at the center of cell faces,  and $\{\mathbf{e}_i\}$ is the standard basis in $\mathbb{R}^d$. We remark that if $x_j$ is away from $\partial \Omega$, then $w_{j}^{i,\pm} = v_{j+1}^{i,\pm}.$ To enforce the zero-flux boundary condition, we set both the velocity $v$ and diffusion $D$ to be zero on $\partial \Omega$. As a result, the columns of $K$ each sum to zero, and the total probability 
$$
\rho^{(l)}\cdot \mathbf{1}=1,\hspace{1cm} \mathbf{1}:= \begin{bmatrix} 1 & \dots & 1
\end{bmatrix}^\top \in \mathbb{R}^N, 
$$ is conserved under time evolution. Since numerical artifacts cause the flux accumulation along the boundary, we also enforce $\rho = 0$ on $\partial \Omega$. When the boundary $\partial \Omega$ is sufficiently far from the trajectory data, this artifact is insignificant. Hereafter, we assume the uniform spatial discretization $\Delta x_i = \Delta x$ for all $i = 1,\dots, d.$ 
Here, we used an explicit time stepping scheme. The Courant--Friedrichs--Lewy (CFL) stability condition enforces $\Delta t = \mathcal{O}(\Delta x^2)$ to ensure the stability of the scheme. To be more specific, we choose
\[
\Delta t < \frac{1}{2d}\frac{\Delta x^2}{D + \Delta x \|v\|_\infty}\,,
\]
where $\|v\|_\infty = \max_i \| v(\mathbf{x})\cdot \mathbf{e}_i\|_\infty$. In this way, we can enforce that all entries of $I + K$ are non-negative with columns summed to one, which implies that $I + K$ is a Markov matrix.

For a complete description of the finite volume scheme, we refer to Ref.~\onlinecite{leveque2002finite}. We remark that there are many higher-order structure-preserving schemes to solve~\eqref{FP2} which also yield a Markov matrix; see Ref.~\onlinecite{HuZhang2022} for example. A more accurate numerical scheme can further reduce the forward modeling error, which is left for future work.

\begin{figure*}
\subfloat[The computed steady state solution to \eqref{eq:FPE} for decreasing values of $\Delta x$.]{
\includegraphics[width = \textwidth]{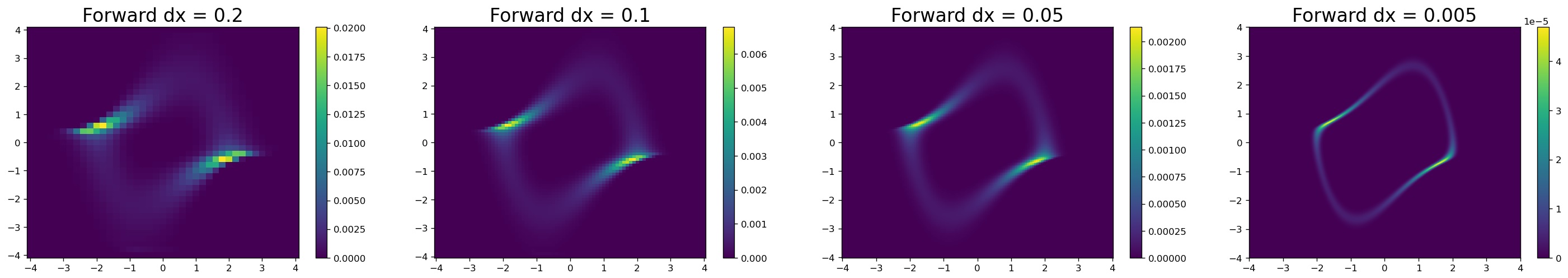}}\\
\subfloat[The approximate physical measure obtained by binning a time trajectory based on the SDE~\eqref{ES} for decreasing values of $\Delta x$.]{\includegraphics[width = \textwidth]{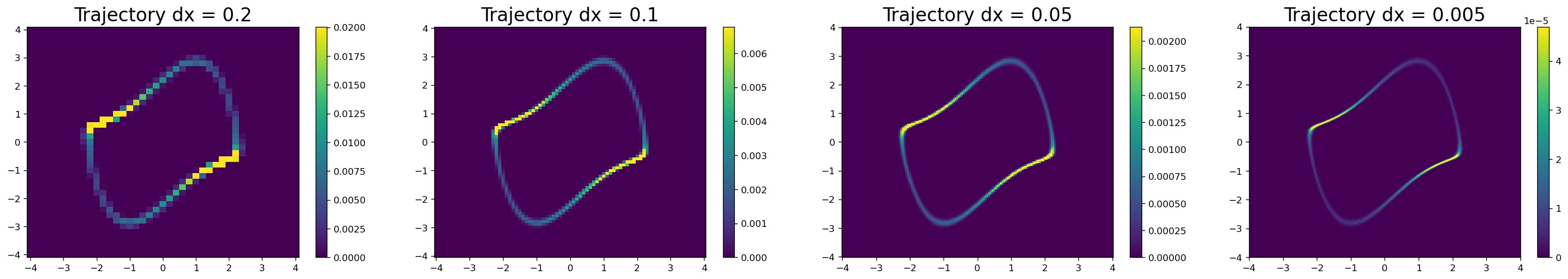}}
\caption{ As the mesh size of the forward model discretization is refined, we visually observe the convergence of the computed steady-state solution (a) to the approximate physical measure (b). The Van der Pol oscillator~\eqref{eq:VDP} with $c = 1$ and $D = 0.001$ is used in this example, and the histograms indicate mass-per cell. \label{F3}}
\end{figure*}

\subsection{Teleportation and Diffusion Regularization}\label{subsec:teleport}
We use the finite volume discretization of the Fokker--Planck equation in~\Cref{sec:FVM} to approximate its steady-state solution. After discretization, finding such stationary distributions to~\eqref{eq:FPE} is equivalent to solving the linear system:
$$
(I+K)\rho =\rho.
$$  
Since the columns of $K$ sum to zero, we have that $M:=I+K$ is a column-stochastic Markov matrix. When $D \neq 0$, $M$ is a transition matrix for an ergodic Markov chain, which has a unique equilibrium. When $D =0$, to guarantee the uniqueness of the equilibrium, Ref.~\onlinecite{yang2021optimal} applies the so-called teleportation regularization~\cite{gleich2015pagerank} and considers
$$M_{\varepsilon}:=(1-\varepsilon) M+ \varepsilon\, U,\quad U = N^{-1} \mathbf{1}\mathbf{1}^{\top}  \in \mathbb{R}^{N\times N}.$$
There is now a unique solution to the linear system 
\begin{equation}\label{eq:tele_ls}
M_{\varepsilon} \rho = \rho, \hspace{.5cm} \rho \cdot \mathbf{1}=1, \hspace{.5cm} \rho >0.
\end{equation}
From a computational aspect, it is useful to take advantage of the fact that $M-I$ is sparse where $I\in\mathbb{R}^{N\times N}$ is the identity matrix and to instead solve
$$(1-\varepsilon)(M-I) \rho  = - N^{-1}\varepsilon \mathbf{1},$$
where we have simply rearranged terms in~\eqref{eq:tele_ls} and used the fact that $\rho \cdot \mathbf{1} = 1$.

Since $U$ is also a column stochastic Markov matrix with the uniform probability of visiting any point of the mesh, using $M_{\varepsilon}$ amounts to stopping the dynamics based on $M$ at a random time and restarting it from a uniformly randomly chosen initial point. The size of $\varepsilon$ represents the restarting frequency--the smaller $\varepsilon$, the rarer we restart.\cite{yang2021optimal}

On the other hand, adding the diffusion component $D$ to the tridiagonal matrix $K$ can be seen as another way of regularizing the noise-free Markov matrix by adding a scaled Brownian motion after each discrete evolution of the deterministic dynamics. For deterministic dynamics with $D=0$, the solution to~\eqref{eq:FPE} might not be unique if there is more than one attractor. The use of teleportation connects all attractors through the ``random restart'', and the solution $\rho_\varepsilon$ to the linear system~\eqref{eq:tele_ls} has support that connects all the disjoint attractors. Similarly, when $D \neq 0$, the Brownian motion connects all disjoint attractors of the deterministic dynamics, giving a unique steady-state solution. In this scenario, the use of teleportation for the diffusive case is simply a numerical treatment to improve the conditioning of matrix $M$ rather than to guarantee the uniqueness of $\rho$.

It is worth noting that both the teleportation regularization and an incorrect diffusion coefficient could be sources of modeling error when we perform parameter identification. Although these regularizations enable faster evaluation of $\rho_\varepsilon(v(\theta))$ and better posedness of the forward problem, they may reduce the accuracy of the inverse problem solution.
\subsection{Numerical Diffusion}\label{sec:numerical_diff}
In Figure~\ref{F3}, we illustrate the $\rho_\varepsilon$ computed as the steady-state solution to the Fokker--Planck equation in the top row and the approximation to physical invariant measures of the corresponding SDE in the bottom row. From Figure \ref{F3}, we see that on a  coarse mesh, the first-order finite volume scheme incurs significant numerical error, which gives a computed solution with an artificial diffusion effect and thus is often referred to as the numerical diffusion.\cite{bewley2012efficient} The amount of numerical diffusion is reduced as the mesh is refined since it is incurred by the first-order scheme. In particular, it is expected to decay as $\pazocal{O}(\max_i{\Delta x_i})$ in the $L^\infty$ norm as we refine the mesh.\cite{leveque2002finite} Besides the teleportation and the modeling diffusion $D$, the presence of numerical diffusion is another modeling error incurred from solving the forward problem.

\section{Gradient Calculation \& Velocity parameterization}\label{sec:gradientcalc}
Another main contribution of this paper is to reconstruct the velocity field $v(\mathbf{x})$ using large-scale parameterizations $v(\mathbf{x};\theta)$, which turns an infinite-dimensional problem of searching for $v(\mathbf{x})$ in a function space to a finite-dimensional optimization problem of finding $\theta \in \Theta \subset \mathbb{R}^m$. Here, we introduce parameterizations based on piecewise-constant, neural network, and global polynomial functions. We also investigate various data-fitting objective functions $\pazocal J$ that compare the mismatch between the observed and simulated invariant measures, $\rho^*$ and $\rho_\varepsilon(v(\theta))$. We compute the gradient of such functions with respect to the coefficients $\theta$ in the parameterized velocity model $v(\mathbf{x};\theta)$ based on the adjoint-state method for the PDE-constrained part and the backpropagation technique~\cite{lecun1989backpropagation} for the neural network part. Thanks to these techniques, we can then efficiently evaluate the gradients of $\pazocal J$ with respect to $\theta$ and thus conveniently use gradient-based optimization algorithms to iteratively update $\theta$, e.g., steepest descent, L-BFGS, conjugate gradient descent methods as well as stochastic methods such as Adam.\cite{kingma2014adam} For notational simplicity, we will write $\rho(v(\theta)) = \rho_{\varepsilon}(v(\theta))$ throughout this section.

\subsection{Gradient Calculation Through the Adjoint-State Method}\label{subsec:adjoint}
Recall the finite volume scheme in~\Cref{sec:FVM} for solving~\eqref{eq:FPE}. The forward model yields a discrete measure $\rho(v(\theta)) = \rho(\theta) = [\rho_1(\theta) \dots \rho_j(\theta) \dots \rho_N(\theta)]^\top$ over the cells $\{C_j\}$, which converges to the solution to~\eqref{eq:FPE} in the weak sense as we refine the discretization parameters. For the explicit form of $\rho(v(\theta))$, we refer to Ref.~\onlinecite[Eqn.~(5.1)]{yang2021optimal}. 
Note that we have highlighted the dependence of our approximate steady-state distributional solution to the Fokker--Planck equation~\eqref{eq:FPE} on the velocity $v(\mathbf{x};\theta)$. Our goal is to solve the optimization problem~\eqref{E2}:
\begin{equation*}
    \inf_{\theta \in \Theta} \pazocal{J}(\rho(v(\theta)),\rho^*)
\end{equation*}
by using gradient-based methods, where $\pazocal{J}$ is the cost function, and $\rho^*$ represents our inference data. The adjoint-state method is an efficient technique by which we can evaluate the derivative $\partial_{\theta} \pazocal{J}$, as the computation time is largely independent of the size of $\theta$. One can derive the adjoint-state method for gradient computations by differentiating the discrete constraint,\cite{nurbekyan2022efficient} which in our case is the eigenvector problem:
$$
g(\rho(\theta),\theta) = M_{\varepsilon}(\theta) \rho(\theta) - \rho(\theta)=\mathbf{0},
$$ 
where $\rho(\theta) \cdot \mathbf{1} = 1.$ Specifically, we compute $\partial_\theta \pazocal{J} = \lambda^\top  \partial_\theta g$ where $ \lambda$ solves $\left( \partial_{\rho} g\right)^\top \lambda =  - \left( \partial_{\rho} \pazocal{J}\right)^\top.$
In our case, this linear system is the adjoint equation (see Ref.~\onlinecite[Eqn.~(5.8)]{yang2021optimal})
\begin{equation}\label{eq:adjoint}
   (M_{\varepsilon}^\top - I)\lambda = - \left( \partial_{\rho} \pazocal{J}\right)^\top +\left( \partial_{\rho} \pazocal{J}\right)^\top \rho \, \mathbf{1}, 
\end{equation}
and the derivative
\begin{equation} \label{eq:adjgradient}
    \partial_\theta \pazocal{J} = \lambda^\top \big(\partial_{\theta} M_{\varepsilon}\big) \rho.
\end{equation}
As a result, we only need to compute the derivatives $\partial_{\rho}\pazocal{J}$ and $\partial_{\theta} M_{\varepsilon}$ to determine the gradient $\nabla_\theta \pazocal{J} = \left(\partial_\theta \pazocal{J} \right)^\top$. The former depends on the choice of the objective function, while the latter is based on a specific parameterization of the velocity field $v(\mathbf{x};\theta)$ determined by its hypothesis space.

\subsubsection{The Computation of $\partial_{\rho}\pazocal{J}$}\label{subsec:objs}
For the objective function $\pazocal{J}$, we consider the quadratic Wasserstein distance, the squared $L^2$ norm, the Kullback--Leibler (KL) Divergence, and the Jensen--Shannon (JS) Divergence. 

\textbf{Quadratic Wasserstein Distance:} For probability measures $\rho$ and $\rho^*$ on $\Omega$, with finite second-order moments, the squared quadratic Wasserstein distance is defined by
    \[
    W^2_2(\rho,\rho^*):=  \inf_{T_{\rho,\rho^*}\in \pazocal{T}}\int_{\Omega}|x-T_{\rho,\rho^*}(x)|^2 d\rho(x),
    \]
    where $$\pazocal{T}:=\{ T:\Omega \to \Omega: \rho(T^{-1}(B)) = \rho^*(B), \, B\in\mathscr{B}\}$$ is the set of maps that push $\rho$ forward into $\rho^*$. \cite{villani2021topics} With an abuse of notation, we also use $\rho(x)$ and $\rho^*(x)$ to denote the densities of $\rho$ and $\rho^*$ respectively. For efficient computation of the $W_2$ distance, we utilize the back-and-forth method,\cite{jacobs2020fast} which instead uses the dual Kantorovich formulation \cite{villani2021topics}
$$W_2^2(\rho,\rho^*)= \sup_{\phi_1,\phi_2}\left( \int_{\Omega} \phi_1(x) \rho^*(x) dx+\int_{\Omega} \phi_2(x) \rho(x) dx
\right),$$
where $\phi_1\in L^1_{\rho^*}(\Omega)$ and $\phi_2 \in L^1_{\rho}(\Omega)$ are required to satisfy  $\phi_1(x)+\phi_2(y) \leq |x-y|^2$. In this case, the Fr\'echet derivative of $\pazocal{J} = W_2^2(\rho,\rho^*)$ with respect to $\rho$ is given by 
$$
\frac{\partial\pazocal{J}}{\partial \rho}= \phi_2.
$$

\textbf{Squared $L^2$ Norm:} The squared $L^2$ distance as the objective function and its Fr\'echet derivative are given by
\begin{align*}
    \pazocal{J} &=\frac{1}{2}\int_{\Omega} |\rho(x) - \rho^*(x)|^2dx, \\
    \frac{\partial \pazocal{J}}{\partial \rho}&=\rho - \rho^*.
\end{align*}
\textbf{KL-Divergence:} The KL-divergence and its Fr\'echet derivative are given by 
\begin{align*}
    \pazocal{J} = D_{\text{KL}}(\rho,\rho^*)&:= \int_{\Omega} \rho^*(x) \log\bigg( \frac{\rho^*(x)}{\rho(x)}\bigg)dx, \\ \frac{\partial D_{\text{KL}}}{\partial \rho} &= -\frac{\rho^*(x)}{\rho(x)}.
\end{align*}

We remark that our definition of the KL-divergence differs from many applications in which it is commonly computed as  $\pazocal{J} = D_{\text{KL}}(\rho^*,\rho).$ 

\textbf{JS-Divergence:} Defining $\rho':=(\rho+\rho^*)/2$, the JS-divergence and its Fr\'echet derivative are given by
\begin{align*}
    \pazocal{J} = D_{\text{JS}}(\rho,\rho^*)& = \frac{1}{2}D_{\text{KL}}(\rho, \rho') + \frac{1}{2}D_{\text{KL}}(\rho^*,\rho'), \\
    \frac{\partial D_{\text{JS}}}{\partial \rho} &= \frac{1}{2}\log\left(\frac{2\rho}{\rho + \rho^*}\right).
\end{align*}
Based on definitions of the KL and JS divergence, it is clear that we may encounter numerical instability issues if either $\rho$ or $\rho^*$ is not supported on the entire domain $\Omega$. Thus, we remark that for the computation of both the KL and JS divergences, we restrict the domain $\Omega$ to regions where both $\rho$ and $\rho^*$ are strictly positive. This is equivalent to the definition of the KL and JS divergence based upon the so-called Csiszar divergence (see Ref.~\onlinecite[Eqn.~(1)]{sejourne2022faster}).

\subsubsection{The Computation of $\partial_{\theta} \pazocal{J}$}\label{subsec:markovdiff}
We have presented a few cases of $\partial_\rho \pazocal J$ for different choices of $\pazocal{J}$. Next, we show how to obtain $\partial_\theta M_\varepsilon$, which is the other necessary component in the adjoint-state method for gradient calculation; see~\eqref{eq:adjoint}-\eqref{eq:adjgradient}. To begin with, we consider $\theta =\{ v_j^i\}$ for all $i=1,\dots,d$ and $j = 1,\dots,N$, which corresponds to one variant of piecewise-constant velocity parameterization.

Since we are only interested in computing the gradient away from $\partial \Omega$, we can utilize the property that $w_j^{i,\pm} = v_{j+1}^{i,\pm}$.
 First, observe that
\[ \frac{\partial M_{\varepsilon}}{\partial v_j^i} = (1-\varepsilon) \sum_{\ell = 1}^d \frac{\Delta t}{\Delta x} \frac{\partial K_{\ell}}{\partial v_j^i} = (1-\varepsilon) \frac{\Delta t}{\Delta x}\frac{\partial K_i}{\partial v_j^i},\]
as well as
\begin{widetext}
\begin{equation}\label{eq:M2}
\frac{\partial K_i}{\partial v_j^i} = \begin{tikzpicture}[baseline={-0.5ex},mymatrixenv]
\matrix [mymatrix,inner sep=4pt](m)
{
\ddots \\& 0 \\ 
\ddots& \vdots\\ &  -H(v^{i}_{j}) & \hspace{.5cm}  -(1-H(v^{i}_{j}))  \\ 
\ddots & \vdots &\vdots \\ 
& H(v^{i}_{j})&\hspace{.5cm} (1-H(v^{i}_{j})) & \hspace{.5cm}0\\
& &  \hspace{.5cm}\vdots & \hspace{.5cm}\vdots & \hspace{.5cm}\ddots\\
& & \hspace{.5cm}0 &\hspace{.5cm}0 \\
& & & \hspace{.5cm}\vdots & \hspace{.5cm}\ddots \\
& & &  \hspace{.5cm}0 \\
& & & & \hspace{.5cm}\ddots \\
};
\mymatrixbraceright{1}{3}{$S_{i}$}
\end{tikzpicture} \hspace{1cm} H(x ):= \begin{cases}
1, & x > 0\vspace{.2cm}\\
0, & x \leq 0
\end{cases}.    
\end{equation}
In \eqref{eq:M2}, $H(\boldsymbol{\cdot})$ is the Heaviside function. We remark that $\partial_{ v_j^i}K_i$ can only be nonzero in the $(j,j)$, $(j,j-S_i)$, $(j-S_i,j)$, and $(j-S_i,j-S_i)$-th entries where $S_i$ is defined in~\eqref{eq:Si}. After solving~\eqref{eq:adjoint} for $\lambda$ and applying~\eqref{eq:adjgradient}, we deduce that
\begin{align}\label{eq:simplifygrad}
\frac{\partial \pazocal{J}}{\partial v_{j}^i} &=  \lambda \cdot \frac{\partial M_{\varepsilon}}{\partial v_j^i}\rho =   (1-\varepsilon)\frac{\Delta t}{\Delta x}\left( \lambda \cdot \frac{\partial K_i}{\partial v_j^i}\rho \right) \nonumber \\ &=(1-\varepsilon)\frac{\Delta t}{\Delta x}\left( H(v^{i}_j)\rho_{j-S_i} \lambda_j+ (1-H(v^{i}_j)) \rho_{j}\lambda_j -H(v_j^{i}) \rho_{j-S_i}\lambda_{j-S_i} - (1-H(v_j^{i}))\rho_j \lambda_{j-S_i} \right)\nonumber \\
 &= (1-\varepsilon)\frac{\Delta t}{\Delta x} \left(\lambda_j - \lambda_{j-S_i}\right)\left(H(v_j^{i})\rho_{j-S_i} +(1-H(v_j^{i})) \rho_{j}\right).
\end{align}
\end{widetext}

Equation~\eqref{eq:simplifygrad} provides an efficient way for computing the gradient of the objective function with respect to the piecewise-constant velocity based on cells  $\{C_j\}$  from our finite-volume discretization.

Alternatively, if the velocity $v = v(\mathbf{x};\theta)$ is smoothly parameterized by the vector $\theta = [\theta_1,\dots,\theta_k,\dots,\theta_m]^\top\in\mathbb{R}^m$, for each $\theta_k$, we can then evaluate
\begin{align}\label{eq:gradient}
    \frac{\partial \pazocal{J}}{\partial \theta_k} &= \sum_{j=1}^N\sum_{i=1}^d \frac{\partial \pazocal{J}}{\partial v_j^i}\frac{\partial v_j^i}{\partial \theta_k}\\  \frac{\partial v_j^i}{\partial \theta_k} &= \mathbf{e}_i\cdot \frac{\partial v}{\partial \theta_k}\bigg\rvert_{(x_{j} - \mathbf{e}_i\Delta x_i/2;\theta)}, \nonumber
\end{align}
to determine the derivative $\partial_{\theta}\pazocal{J}.$ By using a similar indexing convention to~\Cref{sec:FVM}, we can collect the terms $\partial_{v_j^i} \pazocal{J}$ and $\partial_{\theta_k} v_j^i$ into the vectors $\partial_v \pazocal{J}$ and $\partial_{\theta_k} v$, respectively. Therefore, the double summation in \eqref{eq:gradient} is achieved by the inner-product $\partial_v \pazocal{J}\cdot \partial_{\theta_k} v$. Note that for different $\theta_k$, we only need to change $\partial_{\theta_k} v$ as $\partial_v \pazocal{J}$ does not depend on $\theta_k$.

\subsection{Velocity parameterization}\label{sec:velocity_parameter}
We now apply Equations \eqref{eq:simplifygrad} and \eqref{eq:gradient} to evaluate the gradients of several parameterized velocity models. Specifically, we consider piecewise constant, global polynomial, and neural network parameterizations of the velocity.
\subsubsection{Piecewise-Constant parameterization}\label{subsec:PC}
In the case of the piecewise-constant parameterization, we model the velocity as
\begin{equation}\label{eq:PC}
    v(\mathbf{x};\theta)= \sum_{i=1}^d   \sum_{i = j}^N v_j^i \, \chi_{C_j}(\mathbf{x})\, \mathbf{e}_i, \hspace{1cm} \theta = \{v_j^i\}.
\end{equation}
Here, we again use the column-major ordering from~\Cref{sec:FVM} to accumulate vectors of cells $C_j$ with centers $x_j$, and velocity components 
$$
v_j^i=v(x_j - \mathbf{e}_i\Delta x_i /2) \cdot \mathbf{e}_i
$$ 
along the $i$-th direction of the cell face located at $x_j - \mathbf{e}_i\Delta x/2$. The parameter space of the model presented in~\eqref{eq:PC} is given by $\{v_j^i\}$, which has size $N\cdot d$, and the gradient of the parameters $\{v_j^i\}$ can be directly evaluated by~\eqref{eq:simplifygrad}.

We remark that~\eqref{eq:PC} is only one variant of piecewise-constant parameterization since the parameterization mesh is the same as the discretization mesh in the finite-volume method; see~\Cref{sec:FVM}. These two meshes do not have to be coupled together. To reduce the numerical error from the first-order scheme, it is preferable to reduce the spacing $\{\Delta x_i\}$, but we can keep the parameterization mesh fixed so the size of the optimization problem does not change. In this case, we need to apply the chain rule~\eqref{eq:gradient} to obtain the final gradient after evaluating~\eqref{eq:simplifygrad}.

The model defined by~\eqref{eq:PC} can be learned by gradient-based optimization methods. The regularity of the piecewise-constant model defined by~\eqref{eq:PC} can be improved to a $C^0$ function by interpolating between the values $v_j^i$ using either piecewise linear or higher-order piecewise polynomial functions, as in Ref.~\onlinecite{lu2021learning}.

\subsubsection{Global Polynomial parameterization}\label{subsec:GP}
Though the regularity of the piecewise-constant model given by \eqref{eq:PC} can be improved by interpolation, the inverted velocity $v(\mathbf{x};\theta)$ may still be highly oscillatory if the mesh size $\Delta x$ is small. Modeling approaches, such as SINDy, \cite{brunton2016discovering} learn the velocity fields of dynamical systems from a polynomial basis together with sparse regression. Here, we show how the gradient derivation in \eqref{eq:gradient} can be adapted to such polynomial basis parameterization of the velocity field $$v(\mathbf{x};\theta) = [ v^1(\mathbf{x};\theta),\ldots, v^d(\mathbf{x};\theta)]^\top = \sum_{i=1}^d v^i(\mathbf{x};\theta)\, \mathbf{e}_i.$$ The $i$-th component of the velocity field $v^i(\mathbf{x};\theta)$ parameterized by a linear combination of the monomial basis of degree at most $K$ can be written as
\begin{align}\label{eq:poly1}
    v^i(\mathbf{x};\theta) &= \sum_{\ell=1}^{M} a_{\ell}^i (\mathbf{x}^\top \mathbf{e}_1)^{{}^{1}k_{\ell}^i}\dots (\mathbf{x}^\top \mathbf{e}_d)^{{}^{d}k_{\ell}^i}\\ \quad M &= \binom{d+K}{K}, \nonumber
\end{align}
where the powers are represented by multi-indices $$k_\ell^i = ({}^{1}k_{\ell}^i, \dots ,{}^{d}k_{\ell}^i),$$ with $1\leq \ell \leq M$,  $|k_\ell^i| \leq K$, and $\theta = \{a_{\ell}^i \}$. The size of $\theta$ in this case is $d\cdot M$.  To learn the model parameterized by~\eqref{eq:poly1}, we can use \eqref{eq:gradient} to compute the gradient $\partial \pazocal{J} / \partial{a_{\ell}^i}$. Without loss of generality, we assume $\Delta x_i = \Delta x$, for all $1\leq i \leq d$. The only term in~\eqref{eq:gradient} which explicitly depends on the velocity parameterization is
$$
\frac{\partial v_j^i}{\partial a_{\ell}^i} = \left((x_j-\mathbf{e}_i\Delta x/2 )^\top \mathbf{e}_1\right)^{{}^{1}k_{\ell}^i}\dots\left((x_j-\mathbf{e}_i\Delta x/2 )^\top \mathbf{e}_d\right)^{{}^{d}k_{\ell}^i},
$$
where $i$, $j$, $a_\ell^i$ and the multi-index $k_\ell^i$ are fixed. Note that $\partial_{a_{\ell}^{i}} v_j^{i'} = 0$ if $i' \neq i$. Thus, we can again use gradient-based methods to infer proper polynomial coefficients $\{a_{\ell}^i\}$.

Although a global polynomial parameterization guarantees ideal $C^{\infty}$ regularity of the parameterized velocity $v(\mathbf{x};\theta)$, the Runge phenomenon could be a potential downside of this approach. Specifically, as we increase the maximum degree $K$ of the polynomial basis, we may encounter substantial interpolation errors near the boundary $\partial \Omega$.

\subsubsection{Neural Network parameterization}\label{subsec:NNparam}
Motivated by the universal approximation theory of neural networks,\cite{hornik1989multilayer} we may also choose to model each component of the velocity $ v^i(\mathbf{x};\theta)$ as a feed-forward neural network, where the tunable parameters $\theta$ make up the network's weights and biases. We follow Ref.~\onlinecite{li2021bayesian} to combine the adjoint-state method for the PDE constraints and the backpropagation technique to update the weights and biases of the neural network.

The term $\partial_{v_j^i} \pazocal{J}$ in the gradient calculation~\eqref{eq:gradient} can be computed by first evaluating the neural network on the mesh of cell face centers oriented in the direction of $\mathbf{e}_i$ to obtain $\{v_j^i\}$, which is then plugged into~\eqref{eq:simplifygrad} to obtain $\partial_{v_j^i} \pazocal{J}$. The remaining term $\partial_{\theta} v$ in \eqref{eq:gradient} is then computed via the backpropagation technique.\cite{lecun1989backpropagation}

For simplicity, we restrict ourselves to single-layer feed-forward networks. Moreover, by using a smooth activation function, such as the hyperbolic tangent or the sigmoid function, we can guarantee $C^{\infty}$ regularity of the reconstructed velocity $v(\mathbf{x};\theta)$  on the domain $\Omega$. To enforce the zero-flux boundary condition, we manually set $v = 0$ on $\partial \Omega$. Consequently, the neural network parameterization may lack regularity near $\partial \Omega$. However, if the domain is sufficiently large, the support of the physical measure will be very far from $\partial \Omega$, in which case we will not observe any discontinuities originating from the boundary condition while simulating the trajectories based on~\eqref{ES}. As we increase the number of nodes in the hidden layer of the neural network, both the approximation power and the potential difficulty of training the neural network are expected to increase.

\section{Numerical Results}

\begin{figure*}
\subfloat[Ground truth velocity, occupation measure, diffuse trajectory, and non-diffuse trajectory for the Van der Pol oscillator with $c = 0.5$ and $D = 0.02$.\label{fig:VDP_learneda}]{
\includegraphics[width = .95\textwidth]{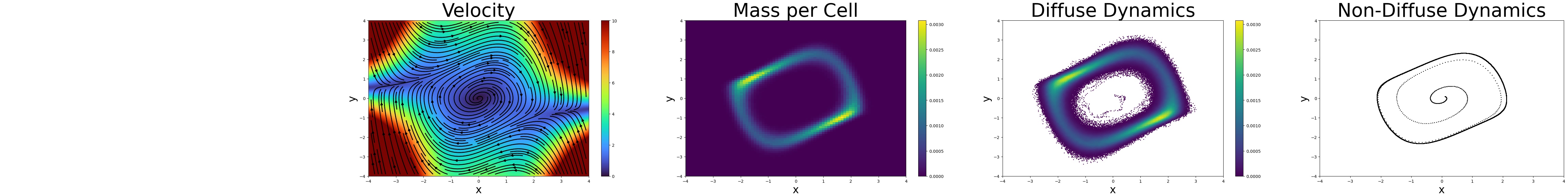}}\\
\subfloat[Piecewise constant parameterization  (see \Cref{subsec:PC}) with the squared $L^2$ objective function.]{
\includegraphics[width = .95\textwidth]{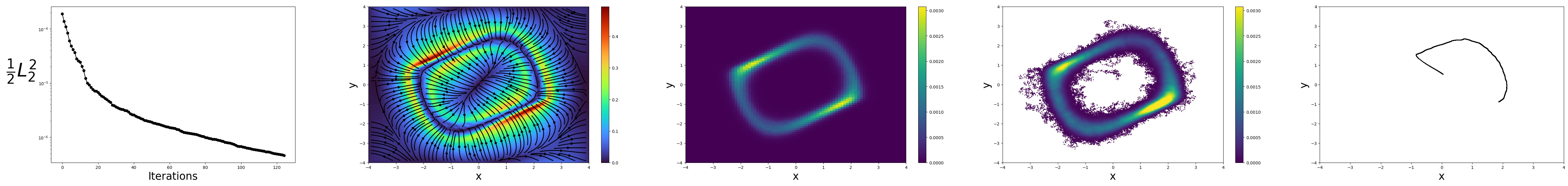}}\\
\subfloat[Degree five global polynomial parameterization  (see \Cref{subsec:GP}) with the squared $L^2$ objective function.\label{fig:VDP_learnedg}]{
\includegraphics[width = .95\textwidth]{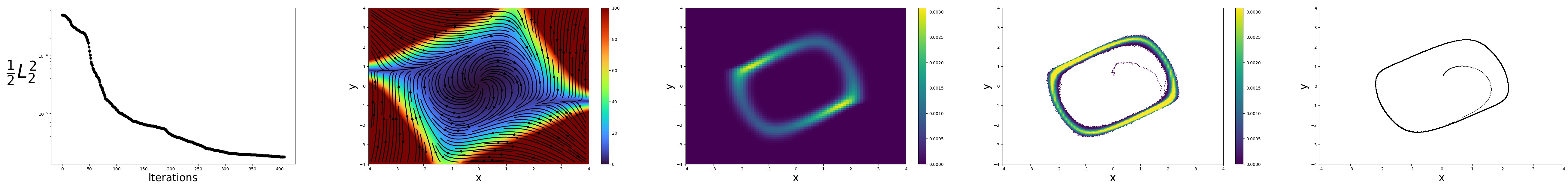}}\\
\subfloat[Neural network parameterization (see \Cref{subsec:NNparam}) with the squared $L^2$ objective function.\label{fig:VDP_learnedb}]{
\includegraphics[width = .95\textwidth]{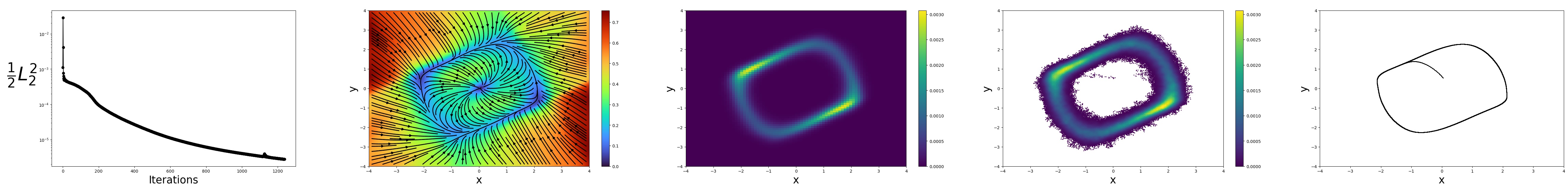}}\\
\subfloat[Neural network parameterization   with the KL divergence objective function.]{
\includegraphics[width = .95\textwidth]{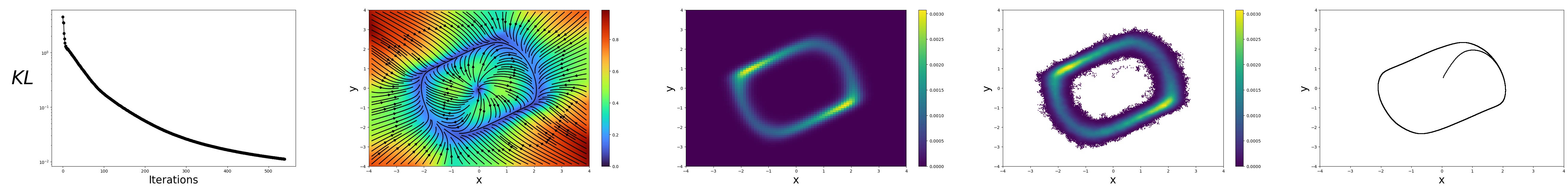}}\\
\subfloat[Neural network parameterization  with the JS divergence objective function.]{
\includegraphics[width = .95\textwidth]{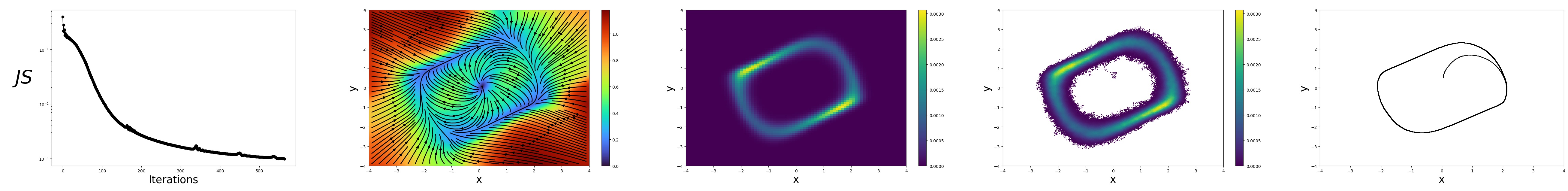}}\\
\subfloat[Neural network parameterization  with the squared $W_2$ objective function.]{
\includegraphics[width = .95\textwidth]{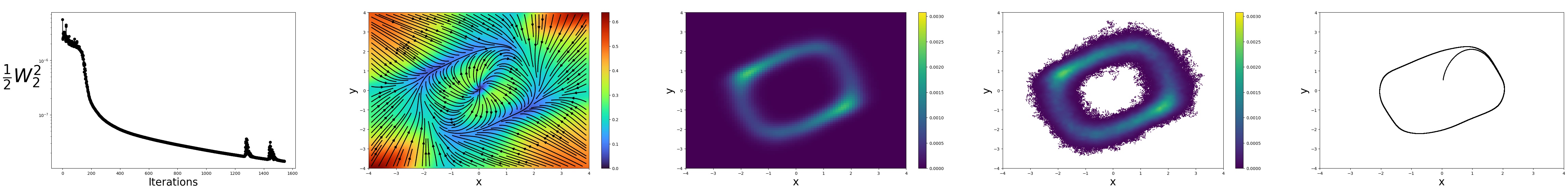}}\\

\caption{Learning velocity fields to reproduce the statistics of the stochastically-forced Van der Pol oscillator. The ground truth occupation measure, velocity, and dynamics are shown in (a). The results for inverting the velocity based on the occupation measure from (a) using neural network, piecewise constant, and global polynomial parameterizations are shown in (b)--(g). The first column shows the objective function; the second column shows the learned velocity vector field; the third column shows the final PDE forward model evaluation based on the learned velocity; the fourth column shows the simulation of a diffuse trajectory, and the final column shows the simulation of a trajectory without diffusion. Specifically, the ``diffuse trajectories'' are simulated according to the Euler--Maruyama method using the assumed diffusion coefficient $D = 0.02$, while the ``non-diffuse'' trajectories assume $D = 0.$ The coloring of each diffuse trajectory is given by the occupation measure it generates; see~\eqref{eq:occupation}. Across all tests, the objective function is minimized to $0.25\%-0.35\%$ of its initial value.
For (b)--(c), the L-BFGS-B algorithm is used for optimization. 
In (d)--(g), the neural network architecture consists of a single hidden layer with the hyperbolic tangent activation function, trained by the Adam optimizer with a learning rate of $10^{-1}$. } \label{fig:VDP_learned}
\end{figure*}

\label{sec:numerics}
In this section, we present several numerical examples to demonstrate the utility of the proposed approach for learning dynamical systems from invariant measures with intrinsic noise.\footnote{We include publicly available code (\href{https://github.com/jrbotvinick/Learning-Dynamics-on-Invariant-Measures}{\texttt{link}}) which contains an example demonstrating the velocity inversion for the Van der Pol oscillator~\eqref{eq:VDP} based on a global polynomial parameterization. It can also be used to reproduce the comparison in Figure \ref{fig:initcompare} and Table \ref{table:compare}.} In \Cref{subsec:vanderpol}, we study the inverse problem for the Van der Pol oscillator with a neural network parameterization of the velocity. In \Cref{subsec:HET}, we time-delay embed a signal sampled from a Hall-effect thruster and proceed to model the dynamics in delay-coordinates based upon the time-delayed invariant measure. We then illustrate that a low-dimensional embedding may increase the uncertainty of the learned model and that the choice of parameterization largely affects the regularity of the reconstructed velocity. In \Cref{subsec:temperature}, we study rolling averages of a temperature data set and perform uncertainty quantification using the learned Fokker--Planck PDE in time-delayed coordinates. We conclude in \Cref{subsec:lorenz} by inverting a component of the Lorenz-63 system's velocity using a neural network parameterization. All experiments are conducted using an Intel i7-1165G7 CPU.

\subsection{Van der Pol Oscillator}\label{subsec:vanderpol}
We begin by considering the autonomous Van der Pol oscillator,\cite{1084738} given by
\begin{equation}\label{eq:VDP}
\begin{cases}
    &\dot{x} = y, \\
    & \dot{y} = c(1-x^2)y - x.
\end{cases}
\end{equation}
Our results for learning a dynamical system with prescribed statistical properties given by the stochastically-forced Van der Pol oscillator are shown in Figure~\ref{fig:VDP_learned}. In the top row, the first panel features the velocity of \eqref{eq:VDP} for the choice of $c = 0.5,$ the second panel shows the approximate occupation measure (see \eqref{eq:occupation}) obtained from the simulation of a single SDE trajectory (see \eqref{ES}), the third panel shows the SDE trajectory used to approximate the invariant measure, and the fourth panel shows the dynamics of the oscillator without stochastic forcing. Throughout, we color the SDE trajectories by their histogrammed density to illustrate the connection between the Lagrangian and Eulerian perspectives. We also stress that the experiment in Figure \ref{fig:VDP_learned} assumes the diffusion coefficient to be known a priori, but that \Cref{subsec:HET} relaxes this assumption.

In the following rows of Figure~\ref{fig:VDP_learned}, we use neural network, piecewise constant, and global polynomial parameterizations of the velocity to solve the inverse problem using the optimization framework from \Cref{sec:FVM,sec:gradientcalc}. For the case of the neural network parameterization, we compare each objective function studied in \Cref{subsec:objs}, while we only focus on the $L^2$ objective for the remaining two parameterizations. Across all tests, the reconstructed velocity is shown to vary significantly from the true velocity shown in the first row of Figure \ref{fig:VDP_learned}. This is mainly due to the lack of data away from the main attracting limit cycle. In regions of the state space with no available data, we can only expect that the modeled velocity $v(\mathbf{x};\theta)$ will direct trajectories towards the attracting limit cycle on which the invariant measure is supported. Indeed, this is what we observe.

\begin{center}
\begin{table}[h!]
\centering

\begin{tabular}{||c| c |c| c||} 
 \hline
 Parameterization &  Objective & Wall-Clock Time (s) &Error  \\ [0.5ex] 
 \hline\hline
  Piecewise constant & $L^2$ & $3.13\cdot 10^1$ & $2.36\cdot 10^{-1}$ \\
 Global polynomial & $L^2$ & $1.75\cdot 10^2$ & $2.90\cdot 10^{-2}$ \\
 Neural network & $L^2$ & $4.14\cdot 10^2$ & $9.07\cdot 10^{-3}$ \\
 Neural network & KL & $2.04\cdot 10^2$ & $7.11 \cdot 10^{-3}$ \\
 Neural network & JS & $2.16\cdot 10^2$ & $9.48 \cdot 10^{-3}$ \\
 Neural network & $W_2$ & $2.16\cdot 10^3$ & $1.07\cdot 10^{-2}$ \\
 \hline
\end{tabular}
\caption{Comparison of the wall-clock computation time and the error for the experiments shown in Figure \ref{fig:VDP_learned}. The error is quantified by the squared $W_2$ distance between the occupation measure of the ground truth diffuse trajectory (see the third panel of Figure \ref{fig:VDP_learneda}), and the occupation measure accumulated from the simulation of a trajectory with diffusion according to the learned velocity (see the fourth panel of Figure~\ref{fig:VDP_learnedb}--\ref{fig:VDP_learnedg}).\label{table:param_quant_compare}}

\end{table}
\end{center}

Moreover, while the learned PDE model \eqref{eq:FPE} matches the observed occupation measure \eqref{eq:occupation} across all tests, we find that the SDE and ODE trajectories generated using the learned velocity vector fields may vary depending on the parameterization. Table~\ref{table:param_quant_compare}
provides a comparison of the accuracy of the learned models, as well as the required computation times. While the piecewise constant velocity is by construction discontinuous and thus does not naturally guarantee the existence and uniqueness of the corresponding ODE solution, the neural network parameterization based on the hyperbolic tangent activation function yields a $C^{\infty}$ velocity. Moreover, while the global polynomial parameterization is also $C^{\infty}$, it may suffer from the Runge phenomena and grow rapidly near the boundary of the domain. Thus, we mainly consider neural network parameterizations of the velocity for the remainder of the numerical tests.

To reduce the computational cost of the inversion in the final row of Figure \ref{fig:VDP_learned}, we compute $\pazocal{J} = W_2^2$ on a coarsened mesh. Among the four objective functions in~\Cref{fig:VDP_learned}, it is worth noting that the $W_2$ metric does not compare the two densities pointwisely and is well-defined for comparing singular measures. The distance reflects both the local intensity differences and the global geometry mismatches.\cite{engquist2020optimal} It has also been shown that the Wasserstein metric is robust to noise.\cite{dunlop2020stability,engquist2020quadratic} Thanks to the geometric nature of the optimal transportation problem, the Wasserstein metric is primarily sensitive to global changes such as translation and dilation and is robust to small local perturbations such as noisy measurements of $\rho^*$. The better stability also brings a downside as the optimization landscape can be relatively flat around the ground truth, which may lead to compromised accuracy in the velocity inversion.

The different velocities shown in the second column of~\Cref{fig:VDP_learned} reveal that there is nonuniqueness if we only use the invariant measure as the reference data. The current modeling assumption yields dynamics reproducing the same invariant measure but does not necessarily recover the same velocity field. Depending on the concrete application, one can add regularization, time information, or focus on velocities in a particular parameterized subspace to avoid nonuniqueness. The large error for the reconstructed velocity near the origin is due to the fact that the method aims to learn the flow on or (in the case of stochastically-forced dynamics) near the invariant measure. It is, therefore, unsurprising that the learned velocity does not match the ground truth where there is no data.

\begin{figure}
    \includegraphics[width = .45\textwidth]{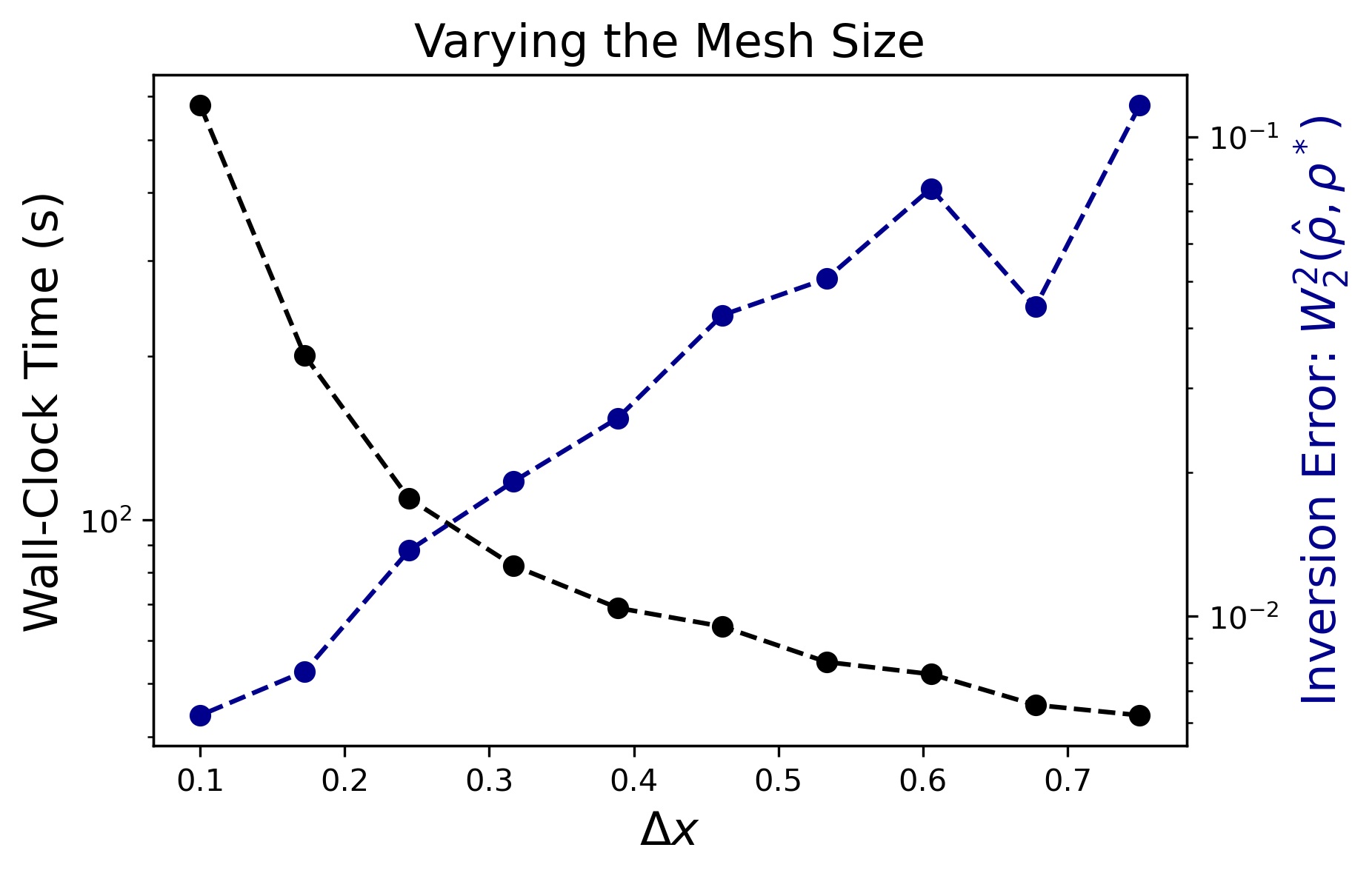}
    \caption{We demonstrate how the computation time and inversion accuracy depend on the mesh spacing used in the first-order FVM solver. Here, we use the Van der Pol oscillator with $D = 0.05$ and learn the velocity using a neural network parameterization. The Adam optimizer is used with a learning rate of $10^{-2}.$ In each case, we reduce the KL divergence objective function below $0.5\%$ of its initial value. The error is quantified in terms of the squared $W_2$ discrepancy between the simulated occupation measure $\hat{\rho}$ according to the learned dynamics and the observed occupation measure $\rho^*$.    \label{fig:meshinv}}

\end{figure}

In Figure~\ref{fig:meshinv}, we show how the inversion accuracy and computation time depend on the chosen value of $\Delta x$. That is, as $\Delta x$ decreases, we can learn velocities that can reproduce the statistics of the observed occupation measure more accurately, with the cost of longer computation time.

Next, we provide experimental details on the comparison of our approach with SINDy~\cite{brunton2016discovering} and the Neural ODE~\cite{DBLP:journals/corr/abs-1806-07366}  frameworks in Figure~\ref{fig:initcompare}. This test uses the Van der Pol oscillator with $c = 2$. Since the SINDy and Neural ODE methods are designed for modeling ODEs, the experiments in Figure~\ref{fig:initcompare} use the diffusion coefficient $D = 0$. While we only plot the first $8$ points of the slowly sampled trajectory in Figure~\ref{fig:initcompare}, the full trajectory used for inference contains $2.5\cdot 10^3$ observations. The quickly sampled trajectory also consists of $2.5\cdot 10^3$ observations. The three approaches considered for comparison have various hyperparameters which can be tuned. For SINDy, we learn the models from the monomial basis up to degree three and use the sequentially thresholded least squares optimizer with threshold $0.025$ to enforce a sparsity condition on the learned coefficients; see Ref.~\onlinecite{brunton2016discovering}. For the Neural ODE framework, the velocity is parameterized by a single-layer fully connected neural network with 100 nodes and a hyperbolic tangent activation function. The Neural ODE is trained using a multiple shooting approach with the mean-squared error objective function. More specifically, rather than treating the simulation of a single long time-trajectory as the forward model, we integrate $N-1$ trajectories initiated at the observed data points $\{\mathbf{x}(t_i)\}_{i=1}^{N-1}$ for a time of $\Delta t = t_{i+1}-t_i$. This approach results in greater success while modeling slowly sampled dynamics. The Adam optimizer with a learning rate of $10^{-3}$ is used, and the tolerance for both relative and absolute error of the ODE solver is set as $10^{-5}$.

To ensure a fair comparison with the Neural ODE framework, we consider our approach based on a neural network parameterization of the velocity using the same architecture, optimizer, and learning rate. For our approach, we use the KL-divergence objective function (see \Cref{subsec:objs}), apply additional Gaussian filtering to the occupation measure (see~\eqref{eq:occupation}) to simplify the resulting optimization, assume a diffusion coefficient of $D = 10^{-3}$ during training, and set $\Delta x = 0.1$. Thus, the only differences between the setup for our approach and the Neural ODE framework are the forward model and objective function.

\begin{center}
\begin{table}[h!]
\centering

\begin{tabular}{||c| c |c| c||} 
 \hline
 Method &  Sampling Freq. & Wall-Clock Time (s) &Error  \\ [0.5ex] 
 \hline\hline
 SINDy & 10.00 & $2 \cdot 10^{-2}$ & $5.6\cdot 10^{-3}$\\
 Neural ODE & $10.00$ & $5 \cdot 10^2 $ & $5.32\cdot 10^{-3} $ \\
 Ours & $10.00$ & $5\cdot 10^2$ & $1.14\cdot 10^{-1} $ \\
 \hline\hline
  SINDy & 0.25 & $10^{-2}$ & $3.52$\\
 Neural ODE & $0.25$ & $5 \cdot 10^2 $ & $1.81$ \\
 Ours & $0.25$ & $5 \cdot 10^2$  &$6.79 \cdot 10^{-2} $ \\
 \hline
\end{tabular}
\caption{Comparison with the SINDy and Neural ODE frameworks for learning from trajectories sampled at different frequencies (Hz). The wall-clock computation time is reported, and the error is quantified by $W^2_2(\hat{\rho},\rho^*)$, where $\hat{\rho}$ is the simulated occupation measure from the learned velocity field, and $\rho^*$ is the observed occupation measure. \label{table:compare}}

\end{table}
\end{center}

As shown in Figure~\ref{fig:initcompare}, all three frameworks can learn from the quickly sampled trajectory. However, SINDy and the Neural ODE frameworks are less robust to changes in the sampling frequency of the inference data than our approach. This is further demonstrated in Table \ref{table:compare}, where we quantify the error in the simulated occupation measure based on the learned velocity. We report the average error over ten trials with different random training seeds to compare our method and the Neural ODE framework. When the data is sampled at a sufficiently high frequency, Table~\ref{table:compare} also shows that methods such as SINDy or the Neural ODE are preferable in terms of both computational cost and accuracy.

\subsection{Hall-Effect Thruster}\label{subsec:HET}

We now turn to the more realistic setting of experimentally sampled time-series data. Specifically, we study the Cathode--Pearson signal sampled from a Hall-effect thruster (HET) in its breathing mode.  Hall-effect thrusters are in-space propulsion devices that exhibit dynamics resembling stable limit cycles while in breathing mode. For details about the experimental setup used to collect the data, the reader is encouraged to consult Refs.~\onlinecite{eckhardt2019spatiotemporal,macdonald2012}. In \Cref{subsubsec:methods}, we utilize Takens' theorem \cite{Takens1981DetectingSA} to reformulate the large-scale optimization framework presented in \Cref{sec:forward,sec:gradientcalc} to be compatible with scalar time-series observations, and in \Cref{subsubsec:results} we demonstrate numerical results based upon this reformulation.

\begin{figure*}
\subfloat[Inference data]{
\includegraphics[width = .325\textwidth]{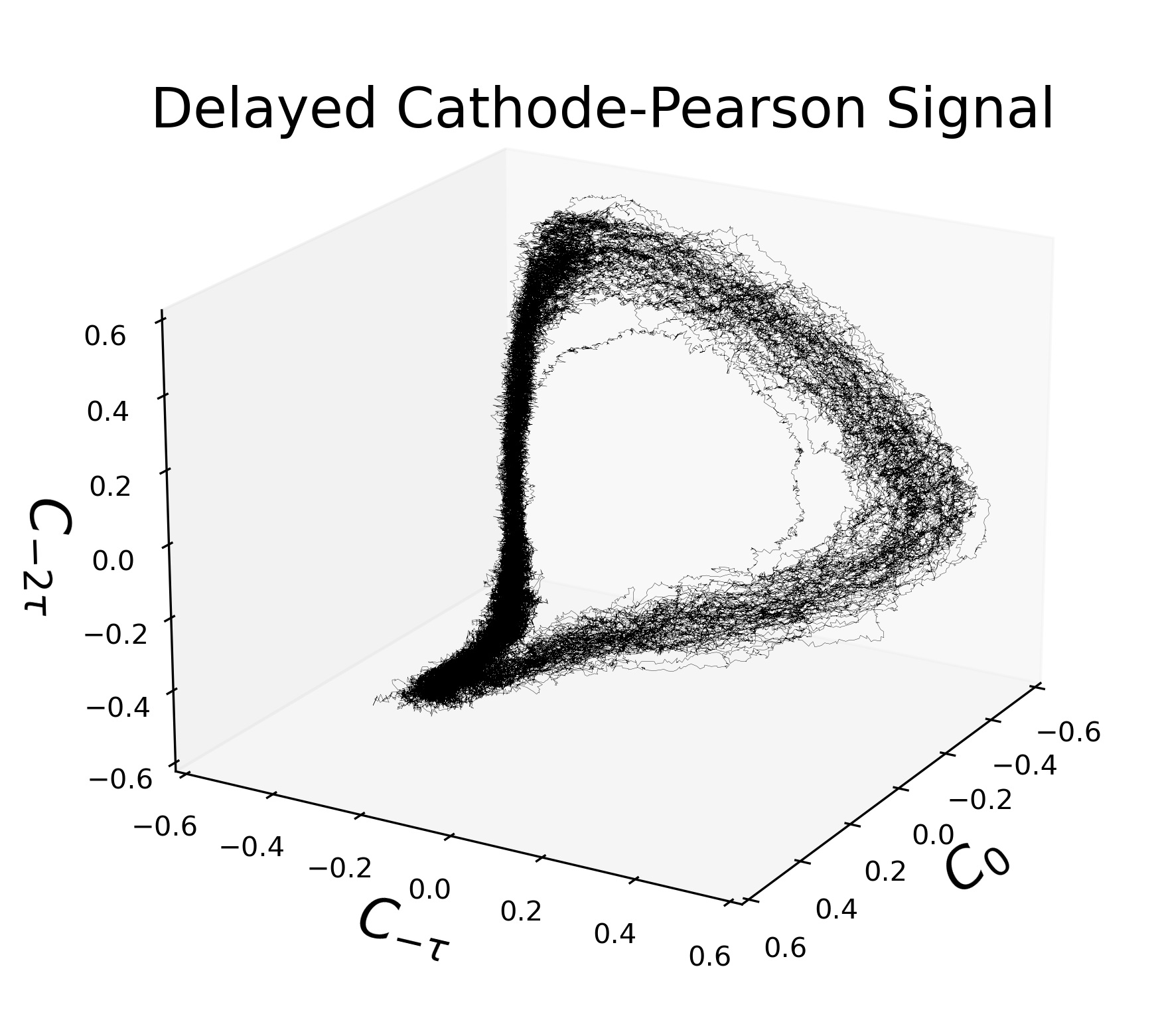}}\label{sf:1}
\subfloat[Learned velocity]{\includegraphics[width = .34\textwidth]{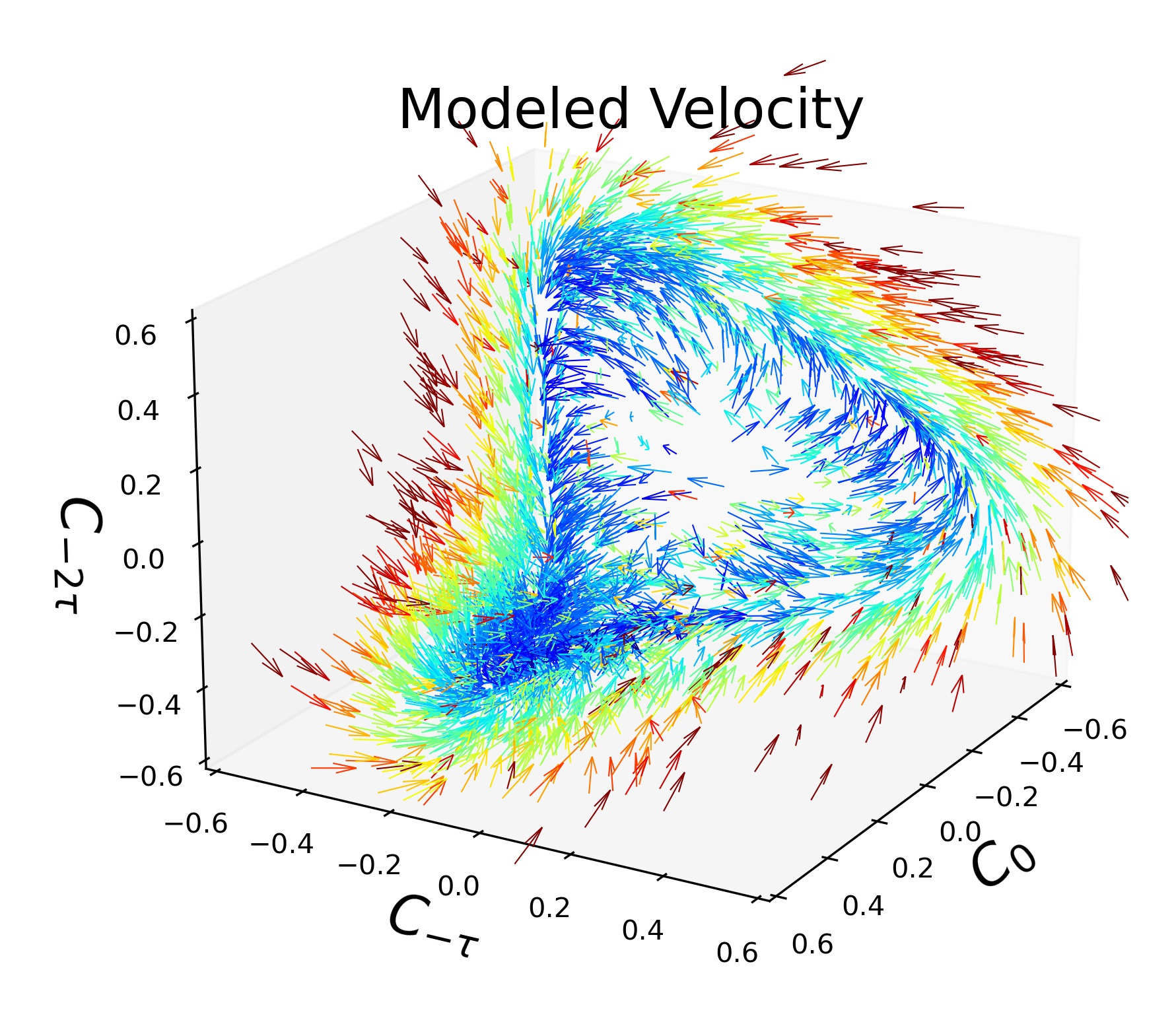}}
\subfloat[Simulated trajectory]{\includegraphics[width = .325\textwidth]{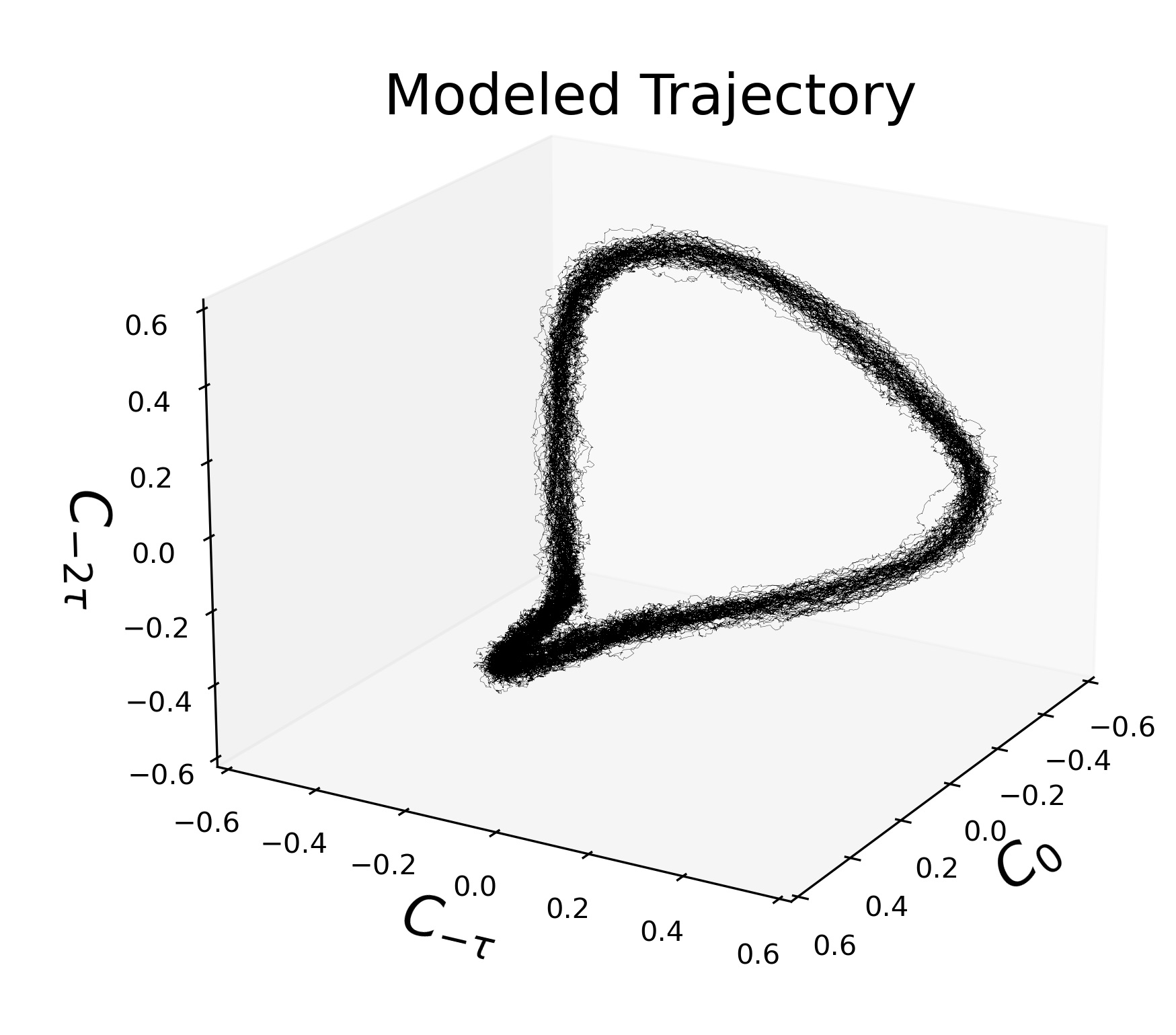}}
\caption{Learning the velocity from the embedded Cathode-Pearson signal's invariant measure. We present the time-delay embedded signal (a), the reconstructed velocity field from the embedded signal's occupation measure (b), and trajectory simulated with the Euler--Maruyama method from the learned velocity and diffusion coefficient (c). In (b), blue indicates slow speed and red indicates fast. The velocity was parameterized by a neural network with $500$ nodes in a single hidden layer and learned using the KL divergence loss function. The three-step procedure in~\Cref{subsubsec:methods} is used to learn the model, and in step one, additional Gaussian filtering is applied to the occupation measure $\rho^*$ to simplify the resulting optimization.\label{fig:HET}}
\end{figure*}

\begin{figure*}
\subfloat[Using the 2D model to predict the evolution of the samples $\mathbf{C}_{2,\tau}$.]{
\includegraphics[width = .95\textwidth]{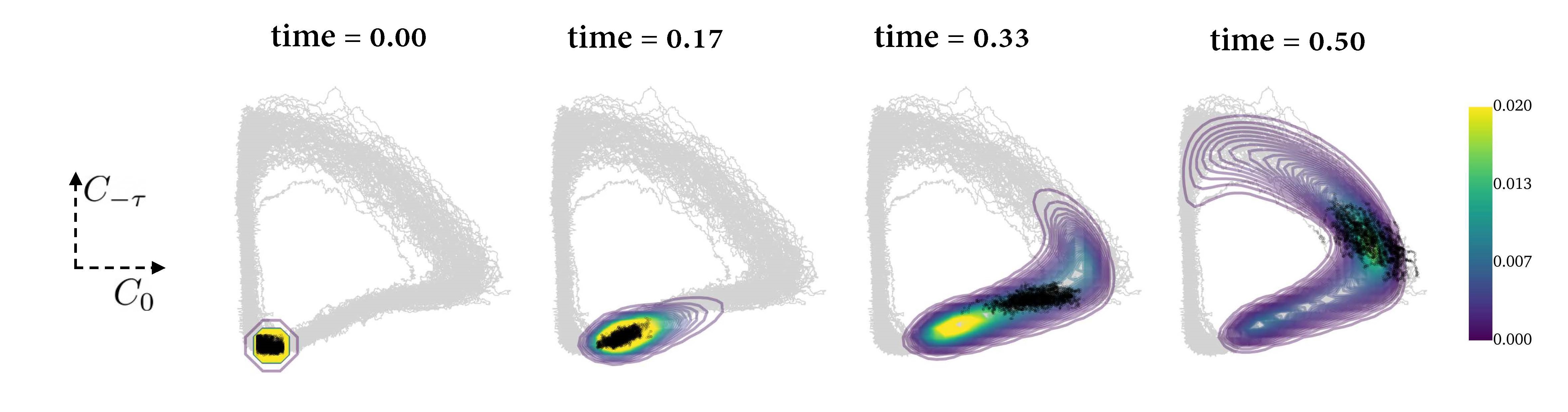}}\\
\subfloat[ Using the 3D model to predict the evolution of the samples $\mathbf{C}_{3,\tau}$.]{\includegraphics[width = .95\textwidth]{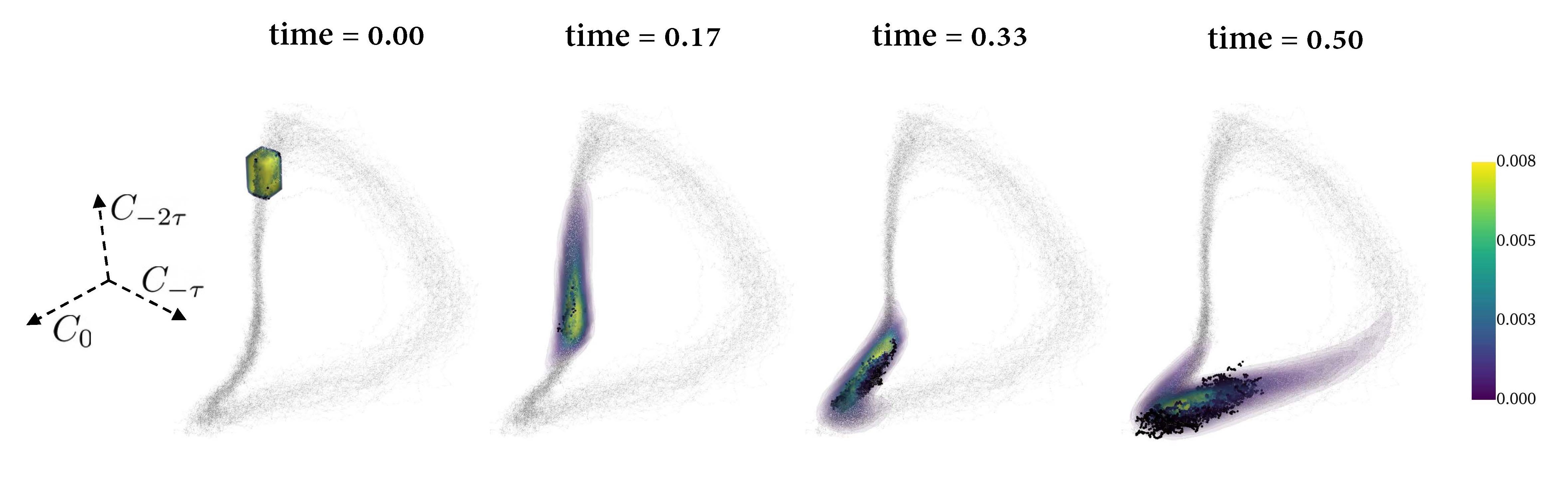}}\\ 
\subfloat[Uncertainty comparison for the 2D and 3D model predictions.]{\includegraphics[width = .95\textwidth]{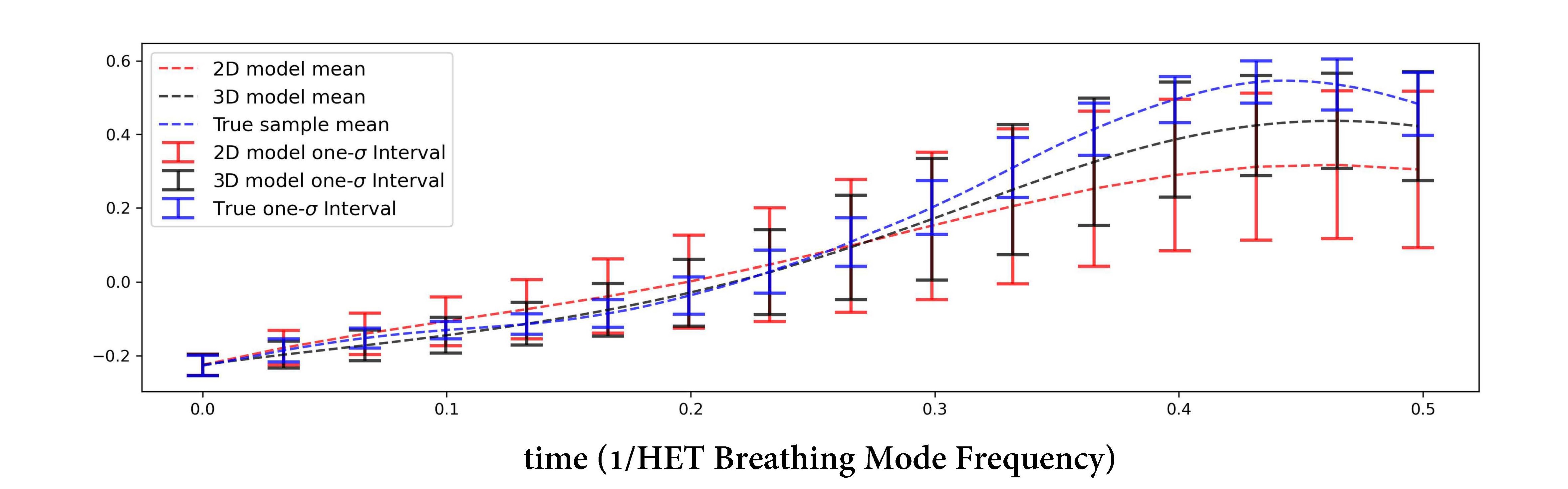}}
\caption{Comparing the model accuracy and uncertainty for the embedded Cathode--Pearson signal with 2D and 3D time delays. The time evolution of the 2D and 3D models is compared to a collection $\{\mathbf{C}_{d,\tau}(t_i)\}_{i=1}^n$ of samples (plotted in black) from the time-delayed Cathode--Pearson signal. The plots (a)--(b) feature a qualitative comparison, whereas (c) shows a quantitative comparison of the uncertainties. Throughout, the time units are normalized to the inverse of a HET breathing mode frequency (16.6kHz). Both the 2D and 3D models utilized a neural network velocity parameterization with 500 nodes in a single hidden layer and reduced the KL divergence objective function to 0.1\% of its initial value during training. As in Figure \ref{fig:HET}, the three-step procedure in~\Cref{subsubsec:methods} is used to learn the models, and in step one, additional Gaussian filtering is applied to the occupation measure $\rho^*$ to simplify the resulting optimization. The 3D visualization was plotted using Ref.~\onlinecite{ramachandran2011mayavi}.}\label{fig:HET_uncertainty}
\end{figure*}

\begin{figure*}
\subfloat[Learned velocity fields for the three parameterizations.]{
\includegraphics[width = \textwidth]{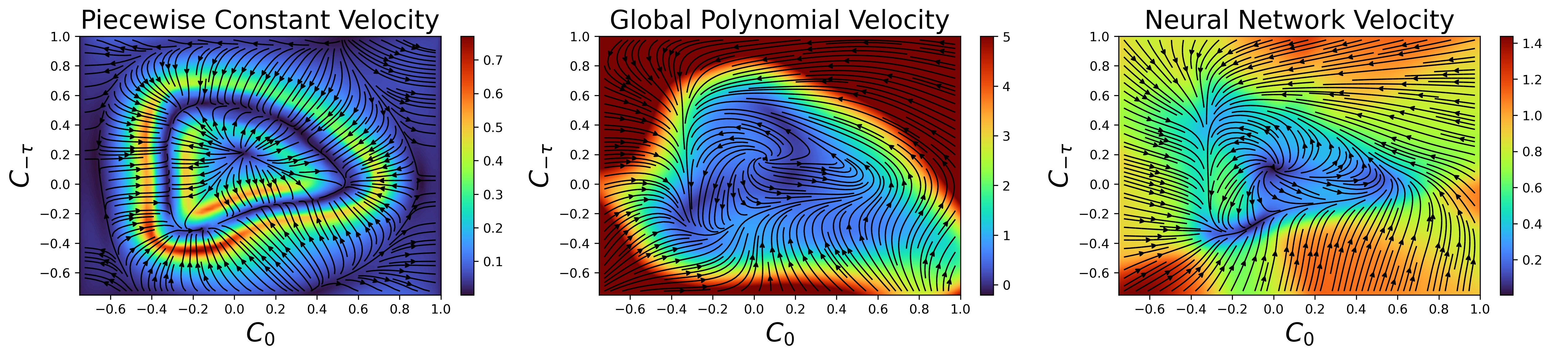}}\\
\subfloat[Close-up view of the learned velocity fields.]{\includegraphics[width = \textwidth]{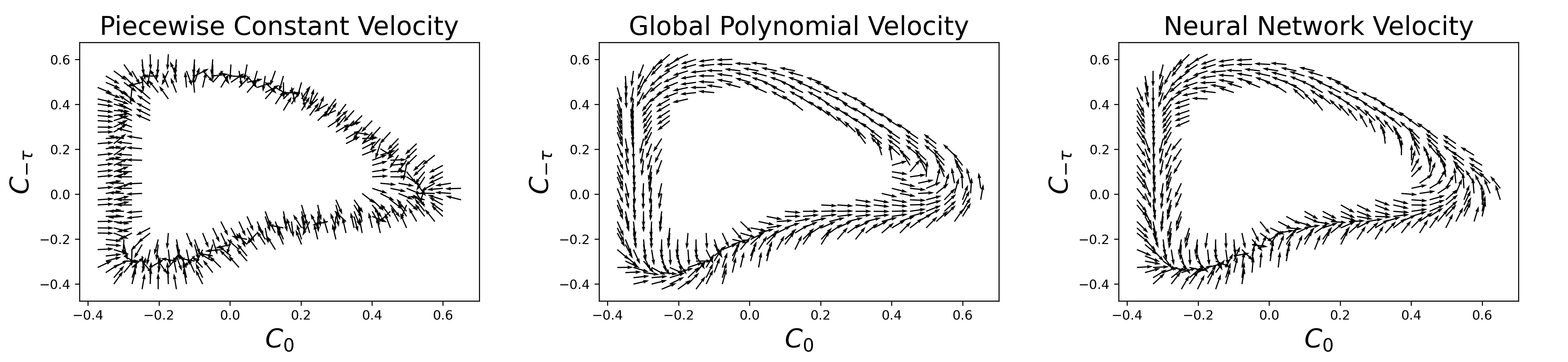}}\\
\subfloat[Forward model evaluations for the learned velocity fields.]{\includegraphics[width = \textwidth]{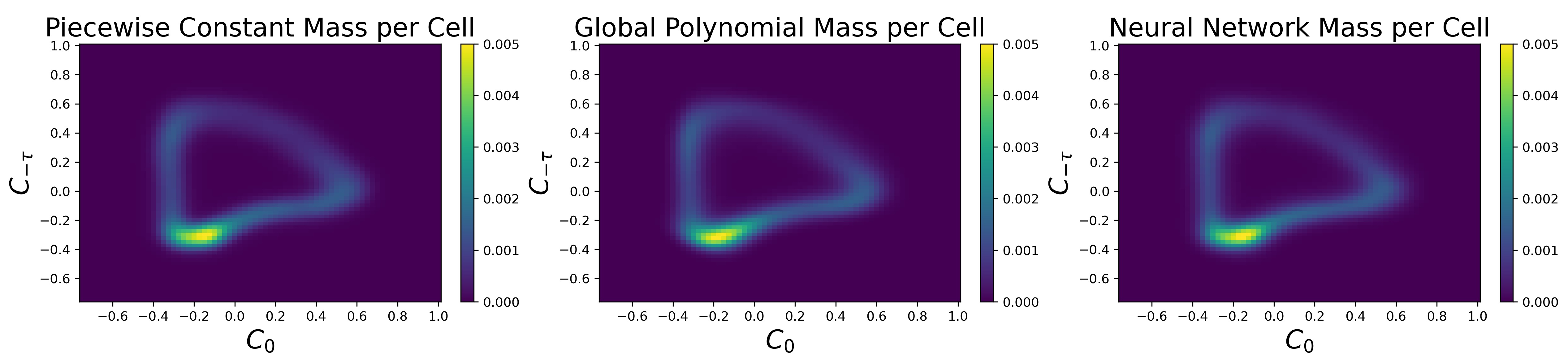}}
\caption{Comparison between the three parameterizations detailed in \Cref{sec:velocity_parameter} for learning a velocity field from the time-delayed Cathode--Pearson signal's invariant measure, using a diffusion coefficient $D = 0.01$. The learned velocities and densities for the piecewise constant (PC), global polynomial (GP), and neural network (NN) discretizations are shown in the three columns, respectively. We show the learned velocity field on the full state space (a), a close-up of the velocity field's direction near the attracting limit cycle (b), and the forward model output $\rho_{\varepsilon}(v(\theta))$ for each parameterization (c). The resulting parameter spaces of these discretizations have a dimensionality of 9800 (PC), 56 (GP), and 400 (NN). The $L^2$ loss is reduced below $0.1\%$ of its initial cost for the PC and NN discretizations and reduced below $0.7\%$ of its initial value for the GP case when we stopped the optimization. \label{fig:comparison}}
\end{figure*}

\begin{figure*}
\subfloat[Occupation measure from weekly rolling averaged temperature measurements.]{
\includegraphics[width = .33333\textwidth]{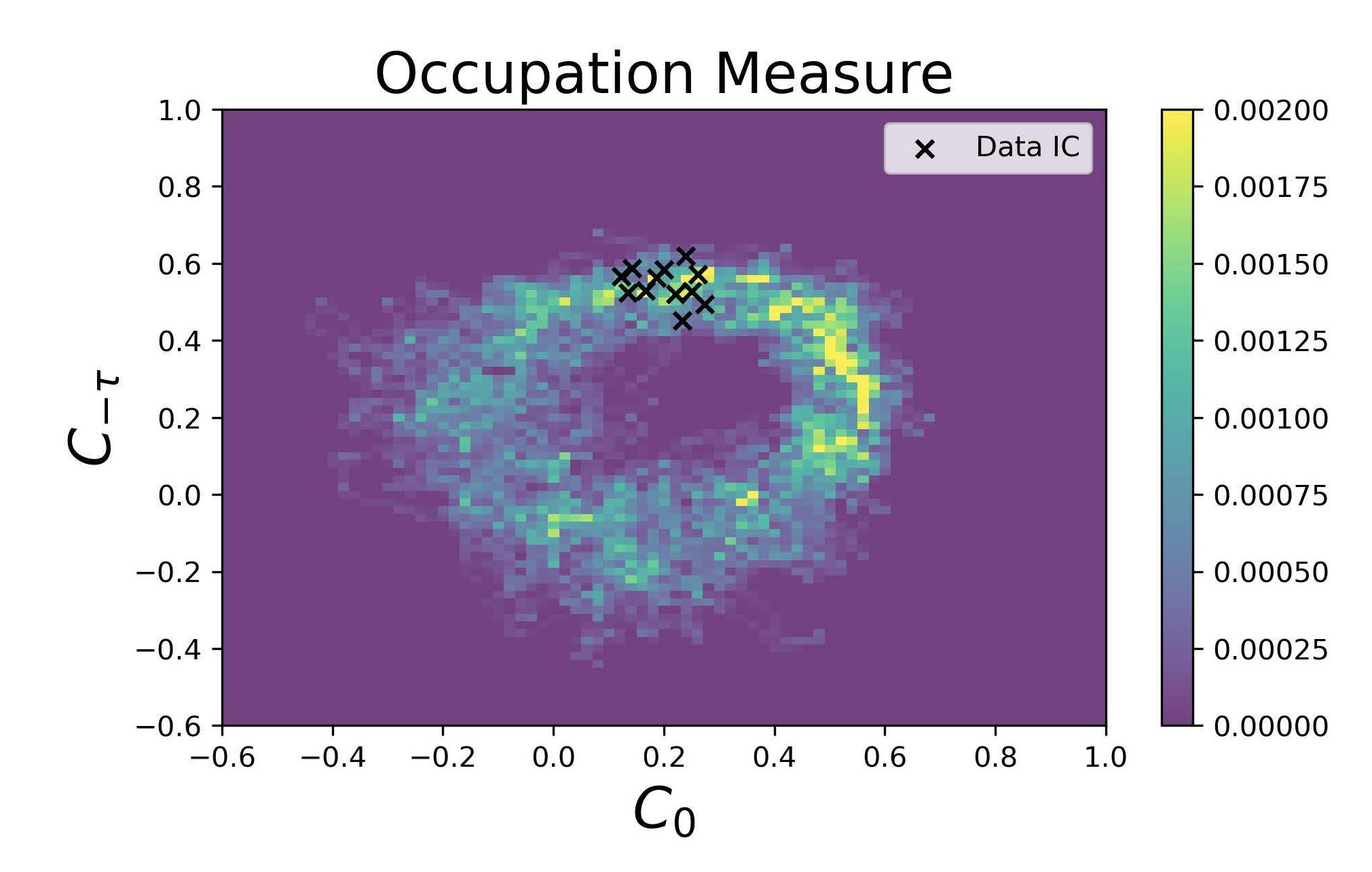}}
\subfloat[Learned velocity from the observed occupation measure.]{\includegraphics[width = .33333\textwidth]{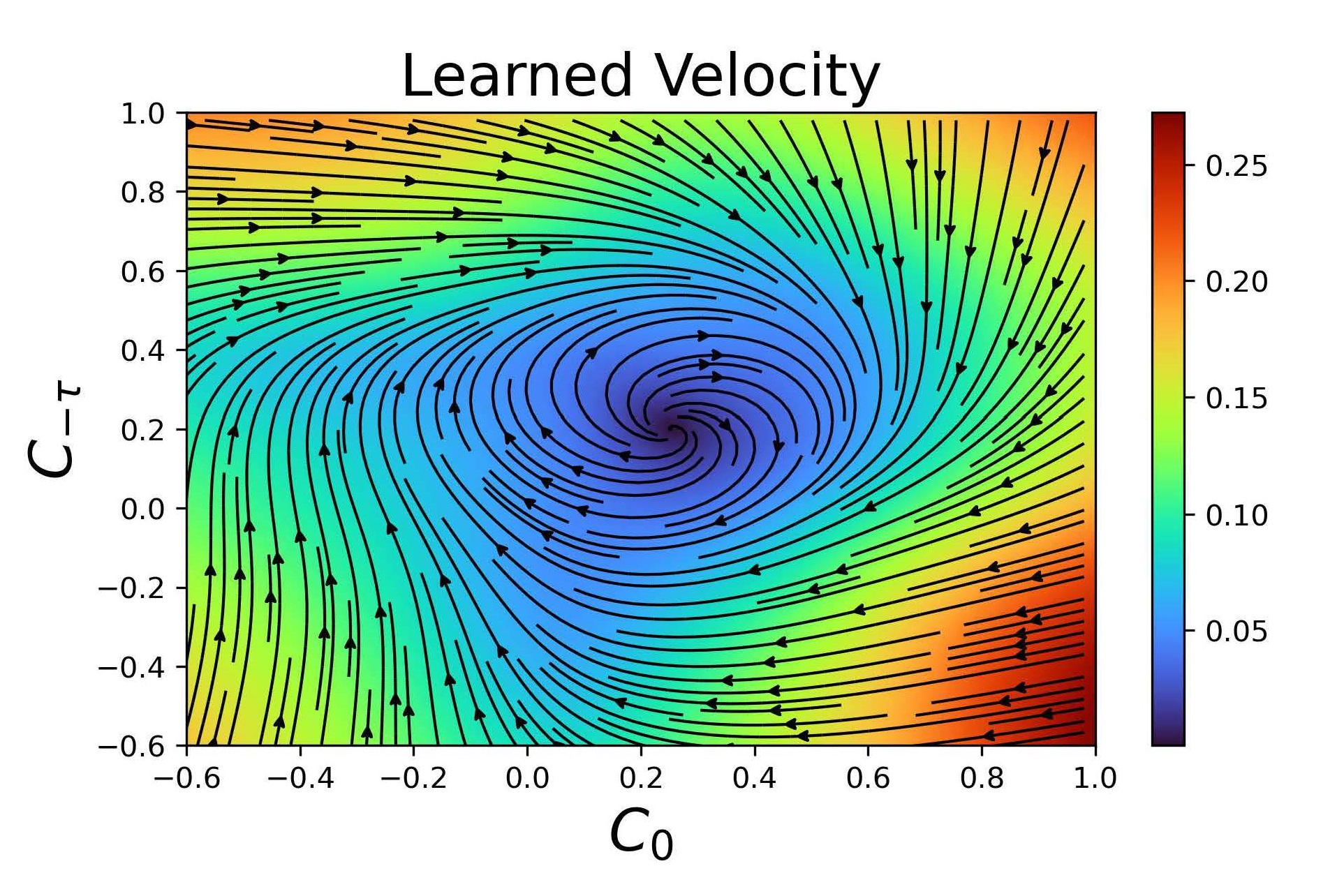}}\subfloat[PDE forward model for the learned velocity ]{\includegraphics[width = .33333\textwidth]{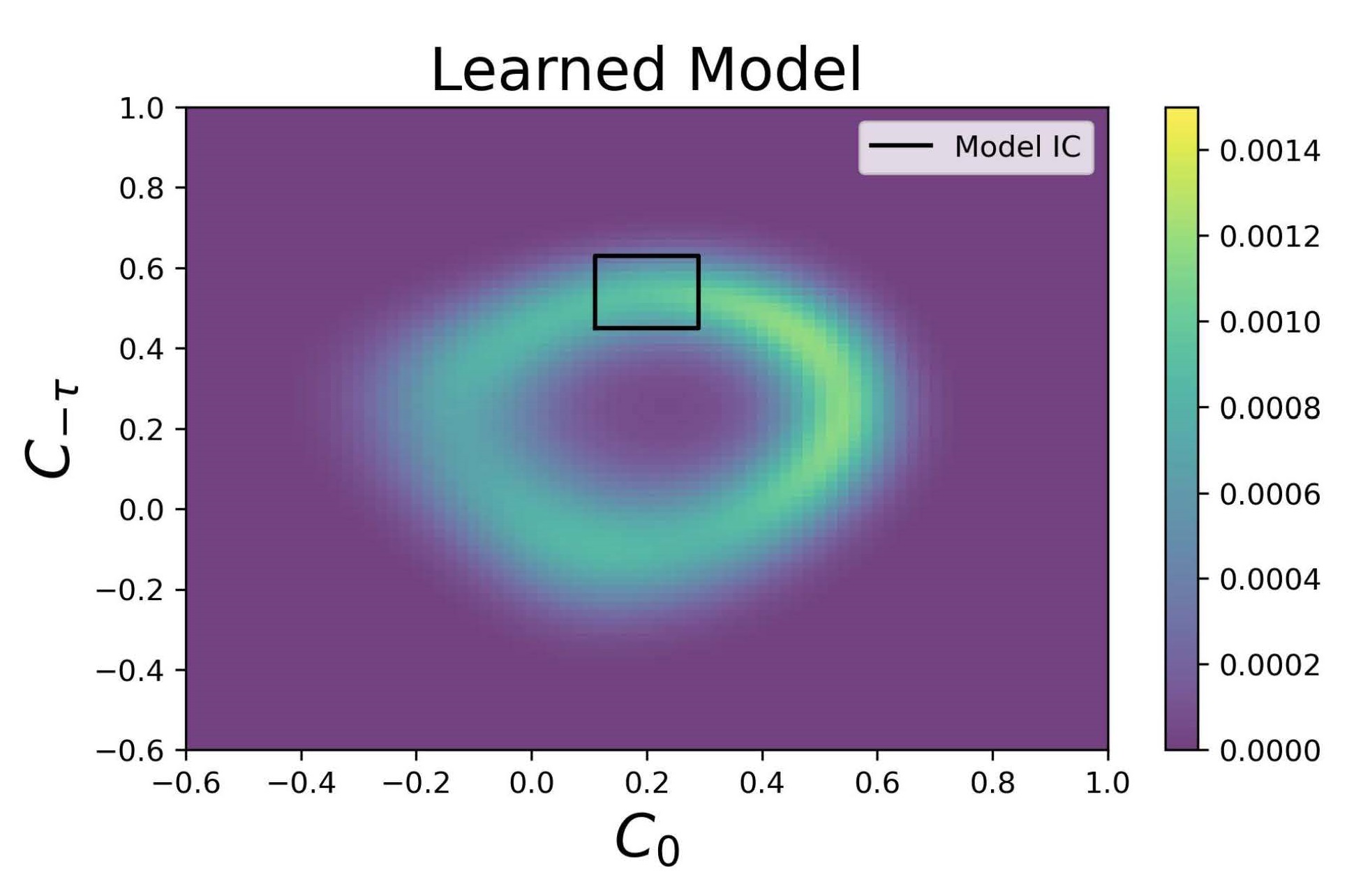}}\\
\subfloat[Evolving the learned Fokker--Planck equation in time to quantify the uncertainty in future temperature forecasts.]{\includegraphics[width = \textwidth]{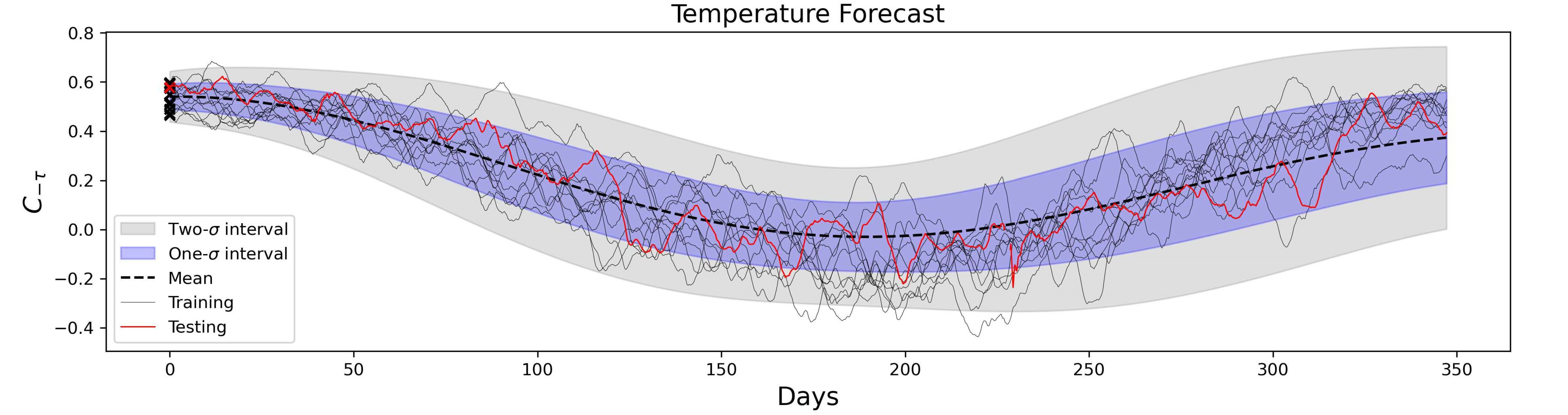}}\\
\caption{Performing prediction and uncertainty quantification for Ithaca, NY's temperature in 2019. (a): The ground truth occupation measure accumulated from 13 years of weekly rolling averaged temperature observations, normalized by an affine transformation to $[-1,1]$; (b): The learned velocity vector field; (c): The corresponding forward model output. In (d), the PDE model with a uniform initialization in the box from (c) is evolved in time and used to quantify the uncertainty in the measurements of $C_0$. Observed trajectories of the temperature in delay coordinates with initial conditions displayed in the top left plot are also shown to demonstrate the effectiveness of the learned model. A time delay of $\tau = 280$ days is used, and the model is trained using a neural network parameterization and the KL-divergence objective function.\label{fig:temp}}
\end{figure*}

\begin{figure*}
\centering
\subfloat[Learned velocity vector field (left), a simulated trajectory with diffusion (middle), and a simulated trajectory without diffusion (right).]{
\includegraphics[width = 1.0\textwidth]{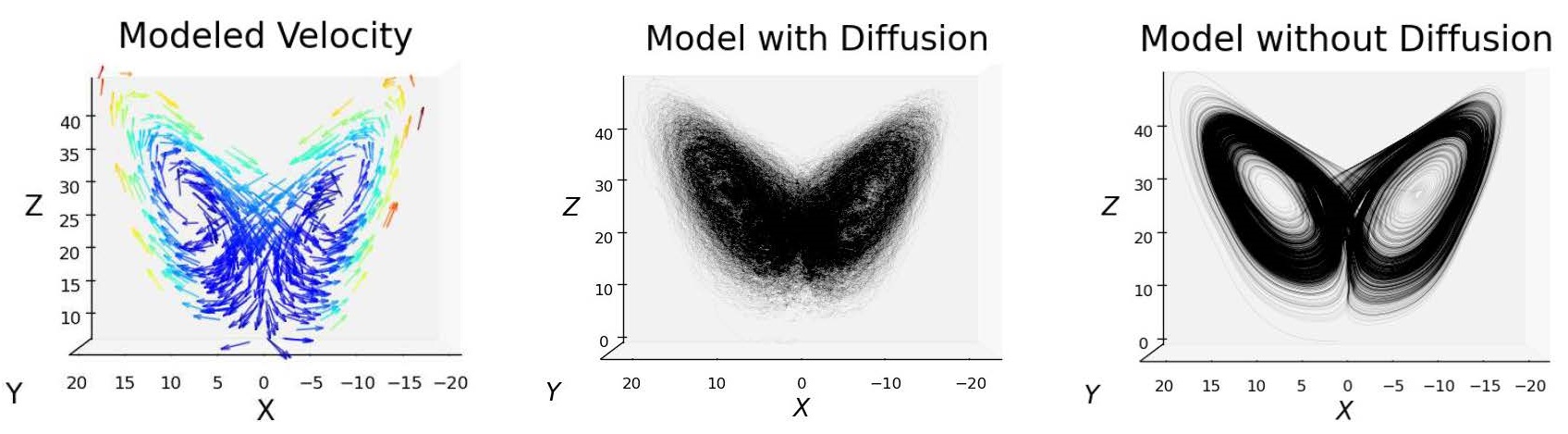}}\\
\subfloat[True velocity (left), a ground truth SDE trajectory (middle), and a ground true ODE trajectory (right). ]{\includegraphics[width = 1.0 \textwidth]{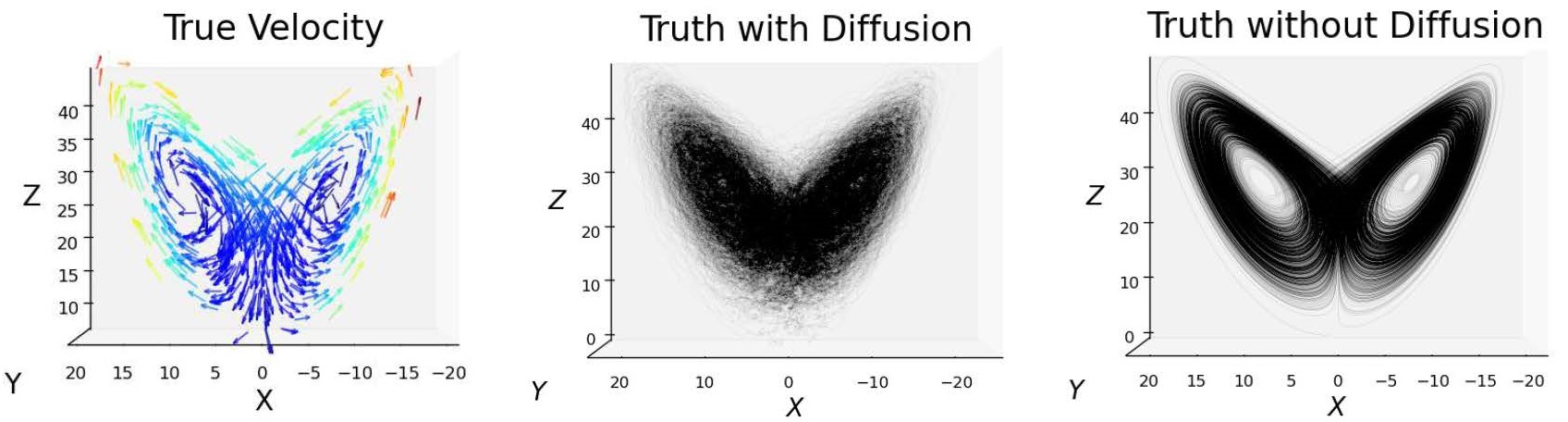}}
\caption{Neural network parameterization of $\dot{x}$ using the Lorenz system's stochastically perturbed invariant measure with $D = 10$ and $\Delta x = 2$. (a): The learned velocity and a single trajectory plotted both with and without diffusion; (b) The ground truth velocity and trajectories both with and without diffusion. For visualization of the occupation measure used to learn the model displayed in the top row, we refer to Ref.~\onlinecite{yang2021optimal}. \label{fig:lorenz}}
\end{figure*}

\begin{figure*}
\centering
\subfloat[Learned velocity vector field (left), a simulated trajectory with diffusion (middle), and a simulated trajectory without diffusion (right).]{
\includegraphics[width = 1.0\textwidth]{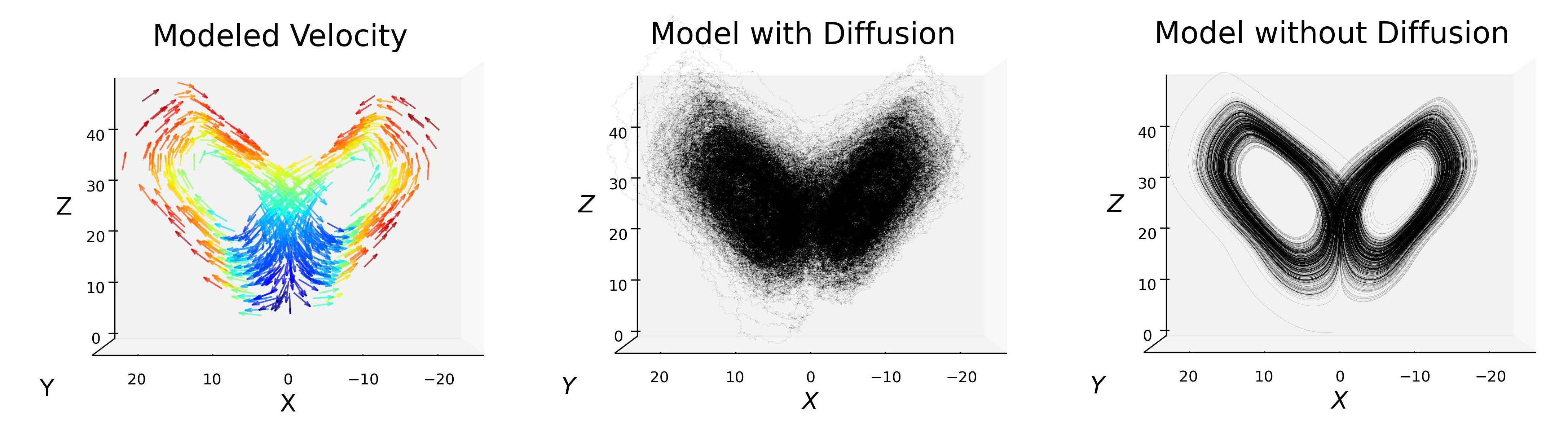}}\\
\subfloat[True velocity (left), a ground truth SDE trajectory (middle), and a ground true ODE trajectory (right).]{\includegraphics[width = 1.0\textwidth]{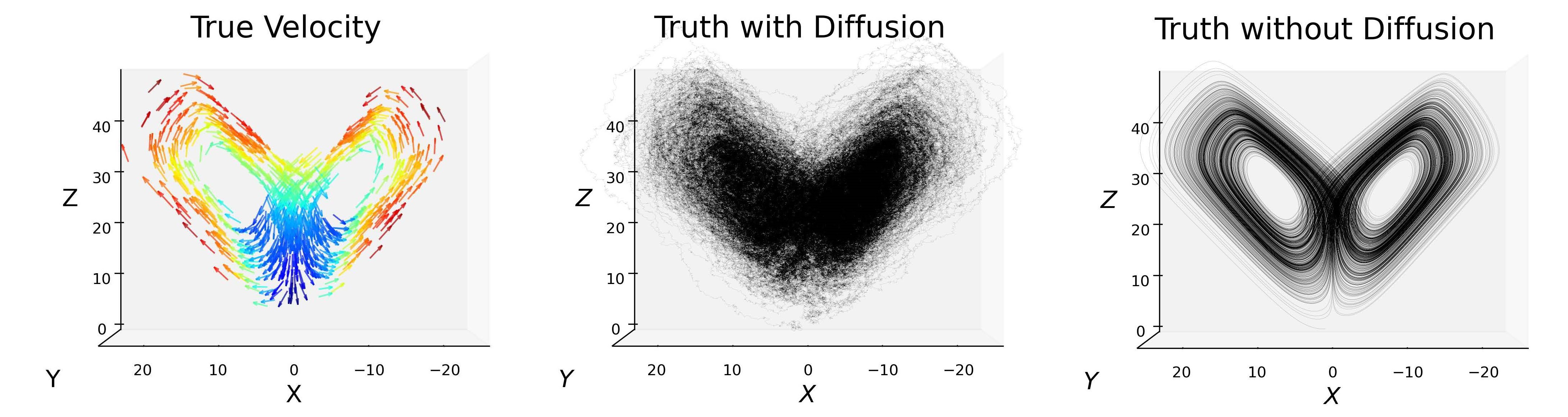}}
\caption{ Neural network parameterization of $\dot{x}$ using the Arctan Lorenz system's stochastically perturbed invariant measure with $D = 10$ and $\Delta x = 2$. The neural network used to learn the velocity contains a single layer of 100 nodes with the sigmoid activation function, and the $L^2$ objective function is used to train the model. In (a), we show results for the learned velocity vector field and the simulation of trajectories using the learned velocity both with and without diffusion. For comparison, in (b), we also include the ground truth velocity and trajectories with and without diffusion. }\label{fig:Alorenz}
\end{figure*}

\subsubsection{Methods}\label{subsubsec:methods}
Intrinsic physical fluctuations present in the Cathode--Pearson signal indicate that the HET's dynamics may be modeled well by a Fokker--Planck equation. Motivated by this insight, we first time-delay embed the Cathode--Pearson signal $C(t)$ in $d$-dimensions to form the trajectory $\mathbf{C}_{d,\tau}(t):=(C(t),C(t-\tau)\dots,C(t-(d-1)\tau))$. We then use a histogram approximation to compute the occupation measure $\rho^*$ of $\mathbf{C}_{d,\tau}(t)$; see \eqref{eq:occupation}. By viewing each dimension of the coordinate system on which the measure $\rho^*$ is supported as the independent variables $C_{-k\tau}(t):=C(t-k\tau)$ where $0\leq k \leq d-1,$ we then seek a solution to the optimization problem \eqref{E2} for a velocity $v = v(\mathbf{C}_{d,\tau};\theta)$. Such a velocity can then provide us with a model of the asymptotic statistics of the embedded trajectory $\mathbf{C}_{d,\tau}(t)$, provided that a suitable diffusion coefficient can be found. 

We note that forming the time-delay coordinates $\mathbf{C}_{d,\tau}(t)$ does require a knowledge of measurements at uniform increments in time. However, the available data may still be sampled slowly enough such that it is impractical to seek a direct approximation of the Lagrangian velocity through the standard approaches described in \Cref{sec:Intro}. This perspective motivates our use of the approach developed in Sections \ref{sec:forward} and \ref{sec:gradientcalc} to learn dynamical systems from invariant measures in time-delay coordinates.

There are a few additional considerations that arise when adapting the modeling framework presented in \Cref{sec:forward,sec:gradientcalc} to real-world data. Namely, we do not know the proper diffusion coefficient a priori (as was the case in \Cref{subsec:vanderpol}). Moreover, the invariant measure that the model is based on does not contain any information about the time scale at which the system evolves. Towards this, we utilize the following three-step procedure as a computationally efficient means to mitigate these difficulties.

\begin{enumerate}
    \item Bin the trajectory $\mathbf{C}_{d,\tau}(t)$ onto a $d$-dimensional mesh with spacing $\Delta x$ along each axis to form the occupation measure $\rho^*$, assume a constant diffusion coefficient $D>0$, and learn the velocity $v = v(\mathbf{C}_{d,\tau};\theta),$ using the framework from \Cref{sec:forward,sec:gradientcalc}. 
    \item Bin the trajectory $\mathbf{C}_{d,\tau}(t)$ onto another $d$-dimensional mesh with spacing $\Delta \hat{x} \leq \Delta x$ to create a new occupation measure $\hat{\rho}^*$ and adjust the diffusion coefficient by solving the optimization problem 
    \begin{equation}\label{eq:newopt}
        \tilde{D} = \argmin_{\hat{D}\in\mathbb{R}} \pazocal{J}(\rho_{\varepsilon}(v;\hat{D}),\hat{\rho}^*),
    \end{equation}
   where the term $\rho_{\varepsilon}(v;\hat{D})$ in \eqref{eq:newopt} denotes the forward model evaluation with the diffusion coefficient $\hat{D}.$  
    \item Rescale both the velocity and diffusion by solving the optimization problem
    \begin{equation}\label{eq:sim}
       \tilde{a}= \argmin_{a \in\mathbb{R}}\sum_{i=1}^N\norm{\hat{\mathbf{C}}(t_i;a) - \mathbf{C}_{d,\tau}(t_i)}_2^2,
    \end{equation}
    where $\hat{\mathbf{C}}(t_i;a)$ denotes the time-$t_i$ solution of the ODE initial value problem with velocity $av(\cdot;\theta)$ and initial condition $\mathbf{C}_{d,\tau}(t_0)$ . The final velocity and diffusion are then given by $\tilde{a} v(\cdot;\theta)$ and $\tilde{a}\tilde{D}$, respectively.
\end{enumerate}

The three-step approach makes repeated use of the fact that $\rho_{\varepsilon}(v;D) = \rho_{\varepsilon}(a v; aD)$, for any scalar multiple $a > 0.$ Indeed, if the true diffusion coefficient $D^* > 0$ is unknown a priori, but we instead seek a solution $v(\cdot,\theta)$ with a different diffusion $ D > 0$, it is guaranteed that the velocity $v = (D/D^*)v^* $ will still provide a solution to the inverse problem. This observation motivates step one, in which an arbitrary diffusion coefficient is used to find a solution $v(\cdot;\theta)$ to the inverse problem. As the dimensionality $d$ is increased, solving the large-scale optimization problem in step 1 on a fine mesh becomes infeasible. As such, step one is typically performed on a coarse mesh where additional Gaussian filtering is applied to the inference measure $\rho^*$ to make the large-scale optimization more feasible.

The diffusion coefficient is then adjusted in step two on a finer mesh via \eqref{eq:newopt} to mitigate the errors due to the Gaussian filtering, numerical diffusion, and histogram errors incurred during step one (see Figure \ref{F3}). Finally, in step three, the scale of both the velocity and diffusion are adjusted via \eqref{eq:sim} such that the time evolution of simulated trajectories is consistent with the inference trajectory $\mathbf{C}_{d,\tau}(t)$ in delay coordinates. Since diffusion plays a relatively small role over short time scales for the quasi-periodic HET data, we use the zero-diffusion trajectory to calibrate a reasonable time-scaling between our model and the available data. However, as the magnitude of the diffusion increases, the least squares fit in \eqref{eq:sim} will become less reliable, and it may be preferable to instead minimize a transport cost between a collection of model samples and a collection of data samples at each time-step. While this final optimization is similar in spirit to various Lagrangian approaches for learning dynamics (see \Cref{sec:Intro}), we remark that the parameter space in \eqref{eq:sim} has only one dimension.

\subsubsection{Results}\label{subsubsec:results}

The results of the three-step procedure in~\Cref{subsubsec:methods} for learning the HET dynamics are shown in Figure \ref{fig:HET} for an embedding dimension of $d = 3$ and time-delay of $\tau = 1.4\cdot 10^{-5}$ seconds, or rather $\tau = .23$ when normalizing the time-scale to the HET breathing mode frequency (16.6kHz). The modeled trajectory accurately reconstructs the shape of the embedded Cathode--Pearson signal but cannot capture the variable diffusion present throughout the time-delayed signal. We do not expect to capture such details, as we assume a constant diffusion coefficient in our model. Nevertheless, we regard the reconstruction of the 3D globally attracting limit cycle as a success and leave extending the model to account for the case of a non-constant diffusion tensor to future work.

The dimensionality of the original HET dynamics is unknown, and as such, a sufficient embedding dimension for the Cathode--Pearson signal is unclear, though likely very high. Interestingly, we can compare the model learned in Figure~\ref{fig:HET} with a 2D analog to demonstrate that when the number of time delays is not sufficiently large, there is more uncertainty in modeling the time-delayed dynamics. This phenomenon is most evident when inspecting regions of the delayed Cathode--Pearson signal for which the 2D embedding lacks structure readily observed in 3D.

Specifically, consider a collection of nearby samples $\{\mathbf{C}_{3,\tau}(t_i)\}_{i=1}^n$ in the 3D time-delay coordinate system $(C_0,C_{-\tau},C_{-2\tau}).$ The corresponding 2D samples $\{\mathbf{C}_{2,\tau}(t_i)\}_{i=1}^n$ will also be nearby one another in the 2D time-delay coordinate system $(C_0,C_{-\tau}).$ In Figure \ref{fig:HET_uncertainty}, we initiate uniform distributions centered about these samples in both 2D and 3D time-delay coordinate systems. We then evolve both the samples and initial uniform distributions forward in time. The evolution of the ground truth samples is simply determined by the time-delayed Cathode Pearson signal $\mathbf{C}_{d,\tau}(t)$, and the evolution of the uniform distributions is given by Fokker--Planck models constructed from the time-delayed Cathode Pearson signal's invariant measure. As the modeled probability densities and ground truth samples evolve in time, we observe in Figure~\ref{fig:HET_uncertainty} that the mean of the 3D model matches the true sample mean more closely than the 2D model and that it has less uncertainty.

In~\Cref{fig:comparison}, we study the three parameterizations from \Cref{sec:velocity_parameter} for learning the time-delayed Cathode--Pearson signal's velocity, now with an embedding dimension of two to allow for clearer visualizations. It can be seen that the density associated with each velocity parameterization indeed matches the ground truth density in~\Cref{fig:comparison}, but that the velocity fields differ significantly from one another. The piecewise-constant velocity in~\Cref{fig:comparison} suffers from poor regularity with discontinuities on the attracting limit cycle. As a result, we lose the connection between the Eulerian and Lagrangian dynamics and cannot reconstruct zero-diffusion trajectories that form a stable limit cycle. On the other hand, the velocities parameterized by the global polynomial and the neural network are both $C^{\infty}$. The differences among these three can clearly be seen via the zoomed-in velocity plots in the second row of~\Cref{fig:comparison}. The global polynomial and neural network discretizations are both global parameterizations of the velocity, and as such, their values near the domain's boundary are dictated by the available data in the center of the domain. This causes the polynomial velocity to rapidly increase near the boundary, and a similar effect can also be seen for the neural network.

It is worth noting that the initial condition for the optimization in~\Cref{fig:comparison} can play a large role in the reconstructed velocity, which is related to the optimization landscape of the nonconvex optimization problem~\eqref{E2} we tackle. In the case of the piecewise-constant discretization, we initialize all velocities to be significantly less than the diffusion coefficient $D = 0.1$. Thus, diffusion initially dominates in the finite volume solver, and all non-boundary cells will contain nonzero mass, which allows for accurate gradient updates everywhere. This phenomenon can also help neural network training, though it is not always necessary due to the global nature of parameterization. Moreover, we initialize our polynomial basis to form the velocity $$(\dot{x},\dot{y}) = (-y+x(0.1-x^2-y^2),x+y(0.1-x^2-y^2)),$$ which describes a globally attracting limit cycle. To converge to the ground truth limit cycle of the time-delayed Cathode--Pearson signal, this initial velocity only needs to be translated and deformed.

\subsection{Temperature Uncertainty Quantification}\label{subsec:temperature}
We now study 2D time-delay embedded data of weekly rolling averages of the temperature in Ithaca, NY, between 2006 and 2020. \cite{USClimateReferenceNetworkafterOneDecadeofOperationsStatusandAssessment} We view temperature fluctuations over short time scales as an intrinsic diffusion process and the approximately periodic oscillation of seasonal temperatures driven by some nonzero velocity. Thus, we model the 2D data in delay coordinates as a diffuse limit cycle. We again follow the procedure in \Cref{subsubsec:methods} to learn a velocity $v(\mathbf{x};\theta)$ and diffusion coefficient $D$, which closely matches the occupation measure.

As in \Cref{subsec:HET}, we can use the trained model $v(\mathbf{x};\theta)$ to quantify measurement uncertainties through the Fokker--Planck equation~\eqref{eq:FPE}, whose solution is a probability density in the time-delay coordinates $(C_0,C_{-\tau})$. Specifically, if we know some initial probability distribution that captures the current state of the temperature system well, we can consider the time evolution of the distribution using our trained model to quantify the uncertainty of future temperature measurements. The process of evolving both the Fokker--Planck PDE from a uniform distribution and the ground truth sample paths from past temperature measurements is shown in Figure~\ref{fig:temp}. The uncertainty bounds from the model accurately capture fluctuations in the training data used to form the occupation measure (plotted in black), as well as a testing sample path previously unseen by the model (plotted in red).

It is also worth noting that the confidence intervals we construct may be larger than the actual range due to several factors, including additional extrinsic noise from filtering the data, modeling errors accumulated from the hypothesis space, numerical diffusion in the forward model, and a sub-optimal embedding dimension. Reducing such errors may result in tighter confidence intervals and considering time delays in higher dimensions could yield better predictions of the temperature's transient behaviors.

\subsection{Lorenz-63 System}\label{subsec:lorenz}
We conclude this section by studying the Lorenz-63 system,\cite{luzzatto2005lorenz} defined by
\begin{equation}\label{eq:lorenz}
\begin{cases}
    \dot{x} &= c_1(y-x) \\
    \dot{y} &= x(c_2-z)-y\\
    \dot{z} &= xy-c_3z
\end{cases},
\end{equation}
where we consider $(c_1,c_2,c_3) = (10,28,8/3).$ For these choices of parameters, the Lorenz-63 system exhibits chaotic behavior and admits a unique physical measure.\cite{TUCKER19991197} In Figure \ref{fig:lorenz}, we assume that the quantities $\dot{y}$ and $\dot{z}$ are known, and we learn a model for the velocity in the $x$-direction, using the stochastically-forced Lorenz-63 system's occupation measure. We emphasize that the data used to approximate the Lorenz system's occupation measure can be sampled slowly or even randomly in time (see Ref.~\onlinecite[Figure 7]{yang2021optimal}). From the approximate occupation measure, we are able to successfully invert the first component $\dot{x}$ of the Lorenz-63 system's velocity via a neural network parameterization. 

We remark that when $\dot{x}, \dot{y},$ and $\dot{z}$ are all simultaneously inverted, the optimization is unsuccessful at reconstructing the true velocity~\eqref{eq:lorenz}. While we may be able to learn a velocity that approximately recovers the stationary state of the Lorenz-63 system in the sense of~\eqref{eq:FPE}, the physical property~\eqref{eq:generic} does not hold. Whether the difficulties of inverting all velocity components of the Lorenz-63 system are due to inherent non-uniqueness in the inverse problem or simply inconvenient local minima during training is worth further investigation in future work. 
To demonstrate the applicability of our approach to non-rational velocities, we also consider the Arctan Lorenz-63 system,\cite{yang2021optimal} given by
\begin{equation}\label{eq:alorenz}
  \begin{cases}
    \dot{x} &= 50\arctan(c_1(y-x)/50) \\
    \dot{y} &= 50\arctan(x(c_2-z)/50-y/50)\\
    \dot{z} &= 50\arctan(xy/50-c_3z/50)
\end{cases}  ,
\end{equation}
where again $(c_1,c_2,c_3) = (10,28,8/3).$ The results for inverting $\dot{x}$ from the occupation measure generated by \eqref{eq:alorenz} with additional stochastic forcing are shown in Figure \ref{fig:Alorenz}, assuming that the quantities $\dot{y}$ and $\dot{z}$ are known.

\section{Conclusion}\label{sec:conclusions}
In this paper, we introduced a PDE-constrained optimization approach to modeling trajectory data originating from stochastic dynamical systems. We first adapted the invariant measure surrogate model in Ref.~\onlinecite{yang2021optimal} based upon the continuity equation to the Fokker--Planck equation. This increased our modeling capacity and prevented overfitting the reconstructed velocity while modeling intrinsically noisy trajectories. We next extended the three-coefficient learning performed in Ref.~\onlinecite{yang2021optimal} to thousands of coefficients by modeling the velocity via global polynomials, piecewise polynomials, and fully connected neural networks. The efficient gradient computation presented in \Cref{sec:gradientcalc} made these large-scale parameterizations of the velocity computationally tractable. We finally studied velocity inversion for invariant measures of time-delay embedded observables. The method of time-delay embedding is useful for analyzing real-world data, where in many cases, only limited observations of complex systems are available. As such, we proceeded to learn the velocity in time-delay coordinates for a Hall-effect thruster system and rolling weekly averages of temperature measurements. Using these models, we predicted future states of the systems and quantified uncertainty in forecasts by evolving the learned Fokker--Planck equation forward in time.

\section*{Acknowledgements}
This paper was supported in part by a fellowship award under contract FA9550-21-F-0003 through the National Defense Science and Engineering Graduate (NDSEG) Fellowship Program, sponsored by the Air Force Research Laboratory (AFRL), the Office of Naval Research (ONR) and the Army Research Office (ARO). R. Martin was partially supported by AFOSR Grants FA9550-20RQCOR098 (PO: Leve) and FA9550-20RQCOR100 (PO: Fahroo). This work was done in part while Y.~Yang was visiting the Simons Institute for the Theory of Computing in Fall 2021. Y.~Yang acknowledges support from Dr.~Max R\"ossler, the Walter Haefner Foundation and the ETH Z\"urich Foundation. This material is based upon work supported by the National Science Foundation under Award Number DMS-1913129. 

We thank Dr.~Chen Li for his helpful suggestions and generosity in sharing code for the approach of \Cref{subsec:NNparam}.

We would like to thank the referees for carefully reading our manuscript and giving many constructive comments that helped improve the paper.

\section*{Author Declarations}
\subsection*{Conflict of Interest}
 The authors have no conflicts to disclose.
\subsection*{Author Contributions}
\textbf{Jonah Botvinick-Greenhouse}: Formal analysis (lead), software (lead), writing -- original draft preparation (equal), writing -- review and editing (equal). \textbf{Robert Martin}: Conceptualization (equal), writing -- review and editing (equal). \textbf{Yunan Yang}: Conceptualization (equal), formal analysis (supporting), supervision (lead), writing -- original draft preparation (equal), writing -- review and editing (equal).
\subsection*{Data Availability}

The data used in Sections \ref{subsec:vanderpol}, \ref{subsec:HET}, and \ref{subsec:lorenz} is available from the authors upon reasonable request. The data used in \Cref{subsec:HET} was obtained from the EPTEMPEST experimental program funded by AFSOR grant FA9550-17QCOR497 (Program Officer: Dr. Brett Pokines). The data used in Section \ref{subsec:temperature} is openly available via Ref.~\onlinecite{USClimateReferenceNetworkafterOneDecadeofOperationsStatusandAssessment}. 

\section*{References}
\bibliography{main.bbl}

\end{document}